# Curvature and symmetry
# of Milnor spheres

By Karsten Grove and Wolfgang Ziller*

*Dedicated to Detlef Gromoll on his $60^{\text{th}}$ birthday*

## Introduction

Since Milnor's discovery of exotic spheres [Mi], one of the most intriguing problems in Riemannian geometry has been whether there are exotic spheres with positive curvature. It is well known that there are exotic spheres that do not even admit metrics with positive scalar curvature [Hi]. On the other hand, there are many examples of exotic spheres with positive Ricci curvature (cf. [Ch1], [He], [Po], and [Na]) and this work recently culminated in [Wr] where it is shown that every exotic sphere that bounds a parallelizable manifold has a metric of positive Ricci curvature. This includes all exotic spheres in dimension 7. So far, however, no example of an exotic sphere with positive sectional curvature has been found. In fact, until now, only one example of an exotic sphere with nonnegative sectional curvature was known, the so-called Gromoll-Meyer sphere [GM] in dimension 7. As one of our main results we prove:

Theorem A. *Ten of the* 14 *exotic spheres in dimension* 7 *admit metrics of nonnegative sectional curvature.*

In this formulation we have used the fact that in the Kervaire-Milnor group, $\mathbb{Z}_{28} = \text{Diff}^+(S^6)/\text{Diff}^+(D^7)$, of oriented diffeomorphism types of homotopy 7-spheres, a change of orientation corresponds to the inverse; hence the numbers 1 to 14 correspond to the distinct diffeomorphism types of exotic 7-spheres.

The exotic spheres that occur in this theorem are exactly those that can be exhibited as 3-sphere bundles over the 4-sphere, the so-called *Milnor spheres*. Each such exotic sphere can be written as an $S^3$ bundle in infinitely many distinct ways; cf. [EK]. Our metrics are submersion metrics on these sphere bundles and we will obtain infinitely many nonisometric metrics on each of these exotic spheres; see Proposition 4.8. We do not know if any of the four remaining exotic spheres in dimension 7 admit metrics with nonnegative curvature, or if the metrics on the Milnor spheres above can be deformed to positive curvature. But no obstructions are known either.

---

*The first named author was supported in part by the Danish National Research Council and both authors were supported by a RIP (Research in Pairs) from the Forschungsinstitut Oberwolfach and by grants from the National Science Foundation.

Another central question in Riemannian geometry is, to what extent the converse to the celebrated Cheeger-Gromoll soul theorem holds [CG]. The soul theorem implies that every complete, noncompact manifold with nonnegative sectional curvature is diffeomorphic to a vector bundle over a compact manifold with nonnegative curvature. The converse is the question which total spaces of vector bundles over compact nonnegatively curved manifolds admit (complete) metrics with nonnegative curvature. In one extreme case, where the base manifold is a flat torus, there are counterexamples [OW], [Ta]. In another extreme case, where the base manifold is a sphere (the original question asked by Cheeger and Gromoll) no counterexamples are known. But there are also very few known examples, all of them coming from vector bundles whose principal bundles are Lie groups or homogeneous spaces (cf. [CG], [GM], [Ri1], and [Ri2]). It is easy to see that the total space of any vector bundle over $S^n$ with $n \leq 3$ admits a complete nonnegatively curved metric. (For $n = 5$, see Proposition 3.14.) Another one of our main results addresses the first nontrivial case.

THEOREM B. *The total space of every vector bundle and every sphere bundle over $S^4$ admits a complete metric of nonnegative sectional curvature.*

The special case of $S^2$ bundles over $S^4$ will give rise to infinitely many nonnegatively curved 6-manifolds with the same homology groups as $\mathbb{C}P^3$, but whose cohomology rings are all distinct; see (3.9).

From a purely topological relationship between bundles with base $S^4$ and $S^7$ (cf. Section 3 and [Ri3]) it will follow that most of the vector bundles and sphere bundles over $S^7$ admit a complete metric of nonnegative curvature; see (3.13). In [GZ2] we will use the constructions of this paper to also analyze bundles with base $\mathbb{C}P^2$, $\mathbb{C}P^2\# \pm \mathbb{C}P^2$, and $S^2 \times S^2$.

From representation theory it is well known that any linear action of the rotation group $SO(3)$ has points whose isotropy group contains $SO(2)$. A proof of this assertion for general smooth actions of $SO(3)$ on spheres was offered in [MS]. However, this turned out to be false. In fact, among other things, Oliver [Ol] was able to construct a smooth $SO(3)$ action on the 8-disc $D^8$, whose restriction to the boundary 7-sphere $S^7$ is almost free, i.e., has only finite isotropy groups. Explicitly, the isotropy groups of the example in [Ol] are equal to $1, \mathbb{Z}_2, D_2, D_3$ and $D_4$. By completely different methods we exhibit infinitely many such actions on the 7-sphere.

THEOREM C (Exotic symmetries of the Hopf fibration). *For each $n \geq 1$ there exists an almost free action of* SO(3) *on $S^7$ which preserves the Hopf fibration $S^7 \to S^4$ and whose only isotropy groups, besides the principal isotropy group* 1, *are the dihedral groups $D_1 = \mathbb{Z}_2$, $D_2 = \mathbb{Z}_2 \times \mathbb{Z}_2$, $D_n$, $D_{n+1}$, $D_{n+2}$, and $D_{n+3}$. Furthermore, these actions do not extend to the disc $D^8$ if $n \geq 4$.*

In the case of the exotic 7-spheres we produce the first such examples.

THEOREM D. *Let $\Sigma^7$ be any (exotic) Milnor sphere. Then there exist infinitely many inequivalent almost free actions of* SO(3) *on $\Sigma^7$, one or more for each fibration of $\Sigma^7$ by 3-spheres, preserving this fibration.*

Since the SO(3) actions in Theorems C and D take fibers to fibers, they induce an action of SO(3) on the base $S^4$. This action of SO(3) on $S^4$ is a fixed action, which yields the well-known decomposition of $S^4$ into isoparametric hypersurfaces [Ca], [HL]. Hence our actions on $S^7$ and $\Sigma^7$ can be viewed as lifts of this action of SO(3) on $S^4$ to the total space of the $S^3$ fibrations.

All of the above results follow from investigations and constructions related to manifolds of cohomogeneity one, i.e., manifolds with group actions whose orbit spaces are 1-dimensional. For closed manifolds this means that the orbit space is either a circle (and all orbits are principal) or an interval. In the first case it is easy to see that the manifold supports an invariant metric with nonnegative curvature. In the second more interesting case, the interior points of the interval correspond to principal orbits and the endpoints to non-principal orbits. Although very difficult, it is tempting to make the following

CONJECTURE. *Any cohomogeneity one manifold supports an invariant metric of nonnegative sectional curvature.*

If true, this would imply in particular that the Kervaire spheres in dimension $4n + 1$, which carry cohomogeneity one actions by SO(2)SO($2n + 1$) (see [HH]), and are exotic spheres if $n$ is even, support invariant metrics of nonnegative curvature. In [BH] it was shown that the singular orbits of the cohomogeneity one actions on the Kervaire spheres have codimension 2 and $2n$, and that they do not admit a metric with positive sectional curvature, invariant under the group action, when $n > 1$.

One of our key results is a small step in the direction of this conjecture.

THEOREM E. *Any cohomogeneity one manifold with codimension two singular orbits admits a nonnegatively curved invariant metric.*

The importance of Theorem E is due to the surprising fact that the class of cohomogeneity one manifolds with singular orbits of codimension two is extremely rich. This is illustrated by our other key result.

THEOREM F. *Every principal $L$ bundle over $S^4$ with $L = $ SO(3) or* SO(4) *supports a cohomogeneity one $L \times$ SO(3) structure with singular orbits of codimension* 2.

Theorem B is now an easy consequence of Theorems E and F in conjunction with the Gray-O'Neill curvature formula for submersions. Theorem A follows from B together with the diffeomorphism classification in [EK]; see the discussion following Remark 4.6. The SO(3) actions in Theorems C and D arise from this construction as well, since the group SO(3) commutes with the principal bundle action and hence induces an action on every associated bundle.

Another consequence of Theorem E is the following (see the discussion after 2.8):

THEOREM G. *On each of the four (oriented) diffeomorphism types homotopy equivalent to $\mathbb{R}P^5$ there exist infinitely many nonisometric metrics with nonnegative sectional curvature.*

The existence of infinitely many cohomogeneity one actions on $S^5$ inducing corresponding actions on each of the homotopy $RP^5$'s was first discovered by E. Calabi (unpublished, cf. [HH, p. 368]), who explained this to us in 1994. These actions can be viewed as the special case $n = 1$ of the Kervaire sphere actions eluded to above. Among these they are the only ones where both singular orbits have codimension 2, so that Theorem E can be evoked directly. Using the same methods as in [Se], one shows that there do not exist any SO(2)SO(3) invariant metrics with positive curvature on these 5-dimensional cohomogeneity one manifolds. This implies that if we apply Hamilton's flow to our metrics of nonnegative curvature, one cannot obtain a metric of positive curvature since Hamilton's flow preserves isometries.

The paper is organized as follows. Section 1 is devoted to general properties of cohomogeneity one manifolds and to an important construction of principal bundles in this framework. In Section 2 we prove Theorems E and G. The constructions in Section 1 are used to prove Theorem F in Section 3. In Section 4 we analyze induced SO(3) actions on associated bundles and derive Theorems C and D. Finally, in Section 5 we examine the geometry of our examples in more detail and raise some open questions.

It is our pleasure to thank J. Shaneson for general help concerning topological questions, and R. Oliver for sharing his insight about SO(3)-actions on discs. We would also like to acknowledge that after seeing a first version of our manuscript in which we had forgotten to include the Calabi examples, H. Rubinstein informed us that O. Dearricott had noticed that Theorem E would yield the existence of metrics of nonnegative curvature on the exotic $RP^5$'s.

## 1. Principal bundles and cohomogeneity one manifolds

We first recall some basic facts about manifolds of cohomogeneity one and establish some notation.

Let $M$ be a closed, connected smooth manifold with a smooth action of a compact Lie group $G$. We say that the action $G \times M \to M$ is of *cohomogeneity one* if the orbit space $M/G$ is 1-dimensional. A cohomogeneity one manifold is a manifold with an action of cohomogeneity one.

Consider the quotient map $\pi : M \to M/G$. When $M/G$ is 1-dimensional, it is either a circle $S^1$, or an interval $I$. In the first case all $G$ orbits are principal and $\pi$ is a bundle map. It then follows from the homotopy sequence of this bundle that the fundamental group $\pi_1(M)$ of $M$ is infinite. In the second case there are precisely two nonprincipal $G$-orbits corresponding to the endpoints of $I$, and $M$ is decomposed as the union of two tubular neighborhoods of the nonprincipal orbits, with common boundary a principal orbit. All of this actually holds in the topological category (cf. [Mo]).

In the remaining part of this paper we will only consider the most interesting case, where $M/G = I$. For this we will make the description above more explicit in terms of an arbitrary but fixed $G$-invariant Riemannian metric on $M$, normalized so that with the induced metric, $M/G = [-1,1]$. Fix a point $x_0 \in \pi^{-1}(0)$ and let $c : [-1,1] \to M$ be the unique minimal geodesic with $c(0) = x_0$ and $\pi \circ c = \mathrm{id}_{[-1,1]}$. Note that $c : \mathbb{R} \to M$ intersects all orbits orthogonally, and $c : [2n-1, 2n+1] \to M$, $n \in \mathbb{Z}$ are minimal geodesics between the two nonprincipal orbits, $B_{\pm} = \pi^{-1}(\pm 1) = G \cdot x_{\pm}$, $x_{\pm} = c(\pm 1)$. Let $K_{\pm} = G_{x_{\pm}}$ be the isotropy groups at $x_{\pm}$ and $H = G_{x_0} = G_{c(t)}$, $-1 < t < 1$, the principal isotropy group. By the slice theorem, we have the following description of the tubular neighborhoods $D(B_-) = \pi^{-1}([-1,0])$ and $D(B_+) = \pi^{-1}([0,1])$ of the nonprincipal orbits $B_{\pm} = G/K_{\pm}$:

$$(1.1) \qquad D(B_{\pm}) = G \times_{K_{\pm}} D^{\ell_{\pm}+1}$$

where $D^{\ell_{\pm}+1}$ is the normal (unit) disk to $B_{\pm}$ at $x_{\pm}$. Hence we have the decomposition

$$(1.2) \qquad M = D(B_-) \cup_E D(B_+)$$

where $E = \pi^{-1}(0) = G \cdot x_0 = G/H$ is canonically identified with the boundaries $\partial D(B_{\pm}) = G \times_{K_{\pm}} S^{\ell_{\pm}}$, via the maps $G \to G \times S^{\ell_{\pm}}$, $g \to (g, \mp\dot{c}(\pm 1))$. Note also that $\partial D^{\ell_{\pm}+1} = S^{\ell_{\pm}} = K_{\pm}/H$. All in all we see that we can recover $M$ from $G$ and the subgroups $H$ and $K_{\pm}$. In fact, two manifolds which carry a cohomogeneity one action by $G$ with the same isotropy groups $H$ and $K_{\pm}$, along a minimal geodesic between nonprincipal orbits, must be $G$-equivariantly diffeomorphic.

In general, suppose $G$ is a compact Lie group and $H \subset K_\pm \subset G$ are closed subgroups such that $K_\pm/H = S^{\ell_\pm}$ are spheres. It is well known (cf. [Bes, p. 95]) that a transitive action of a compact Lie group $K$ on a sphere $S^\ell$ is linear and is determined by its isotropy group $H \subset K$. Thus the diagram of inclusions

$$(1.3) \qquad \begin{array}{ccccc}
 & & G & & \\
B_- = G/K_- & \overset{j_-}{\diagup} & & \overset{j_+}{\diagdown} & B_+ = G/K_+ \\
 & K_- & & K_+ & \\
S^{\ell_-} = K_-/H & \underset{i_-}{\diagdown} & & \underset{i_+}{\diagup} & S^{\ell_+} = K_+/H \\
 & & H & &
\end{array}$$

determines a manifold

$$(1.4) \qquad M = G \times_{K_-} D^{\ell_-+1} \cup_{G/H} G \times_{K_+} D^{\ell_++1}$$

on which $G$ acts by cohomogeneity one via the standard $G$ action on $G \times_{K_\pm} D^{\ell_\pm+1}$ in the first coordinate. Thus the diagram (1.3) defines a cohomogeneity one manifold, and we will refer to it as a cohomogeneity one group diagram, which we sometimes denote by $H \subset \{K_-, K_+\} \subset G$. We also denote the common homomorphism $j_+ \circ i_+ = j_- \circ i_-$ by $j_0 \colon H \to G$.

We are now ready for the main construction in this section: *Principal bundles over cohomogeneity one manifolds.*

Let $L$ be any compact Lie group, and $M$ any cohomogeneity one manifold with group diagram $H \subset \{K_-, K_+\} \subset G$. It is important to allow the $G$-action on $M$ to be noneffective, i.e. $G$ and $H$ have a common normal subgroup, since this will produce more principal bundles over $M$; see, e.g., (3.1), and (3.2).

For any Lie group homomorphisms $\phi_\pm : K_\pm \to L$, $\phi_0 : H \to L$ with $\phi_+ \circ i_+ = \phi_- \circ i_- = \phi_0$, let $P$ be the cohomogeneity one $L \times G$-manifold with diagram

$$(1.5) \qquad \begin{array}{ccccc}
 & & L \times G & & \\
 & \overset{(\phi_-,j_-)}{\diagup} & & \overset{(\phi_+,j_+)}{\diagdown} & \\
K_- & & & & K_+. \\
 & \underset{i_-}{\diagdown} & & \underset{i_+}{\diagup} & \\
 & & H & &
\end{array}$$

Clearly the subaction of $L \times G$ by $L = L \times \{e\}$ on $P$ is free since

$$L \cap (l,g)K_\pm(l,g)^{-1} = (l,g)(L \times \{e\} \cap K_\pm)(l,g)^{-1}$$

as well as

$$L \cap (l,g)H(l,g)^{-1} = (l,g)(L \times \{e\} \cap H)(l,g)^{-1}$$

is the trivial group for all $(l, g) \in L \times G$. Moreover, $P/L = M$ since it has a cohomogeneity one description $H \subset \{K_-, K_+\} \subset G$. It is also apparent that the nonprincipal orbits in $P$ have the same codimension as the nonprincipal orbits in $M$, as well as the same slice representation, since the normal bundles in $M$ pull back to the normal bundles in $P$ under the principal bundle projection $P \to M$. In summary:

PROPOSITION 1.6. *For every cohomogeneity one manifold $M$ as in (1.3) and every choice of homomorphisms $\phi_\pm : K_\pm \to L$ with $\phi_+ \circ i_+ = \phi_- \circ i_-$, the diagram (1.5) defines a principal $L$ bundle over $M$.*

Note, moreover, that the $L \times G$-action on $P$ may well be effective even if the $G$-action on $M$ is not.

We now move on to discuss *induced actions on associated bundles*:

Let $F$ be a smooth manifold on which $L$ acts, $L \times F \to F$. Consider the total space of the associated bundle $V = P \times_L F$. Observe that the product of the trivial $G$-action on $F$ with the sub-action of $G = \{e\} \times G \subset L \times G$ on $P$ induces a natural $G$-action on $V$.

LEMMA 1.7 (Isotropy Lemma). *The natural $G$-action on $V = P \times_L F$ has exactly the following types of isotropy groups*

$$\phi_\pm^{-1}(L_u) \quad and \quad \phi_0^{-1}(L_u)$$

*where $L_u$, $u \in F$ are the isotropy groups of $L \times F \to F$.*

*Proof.* Consider the $L$-orbit, $L(x, u) = \{(\ell x, \ell u) \mid \ell \in L\}$ of a point $(x, u) \in P \times F$. Then

$$
\begin{aligned}
G_{L(x,u)} &= \{g \in G \mid gL(x, u) = L(gx, u) = L(x, u)\} \\
&= \{g \in G \mid \exists\, \ell \in L : (gx, u) = (\ell x, \ell u)\} \\
&= \{g \in G \mid \exists\, \ell \in L_u : (\ell^{-1}, g) \in (L \times G)_x\}.
\end{aligned}
$$

However, $(L \times G)_x$ is some $(\hat\ell, \hat g)$-conjugate of one of $(\phi_+, j_+)(K_+), (\phi_-, j_-)(K_-)$ or $(\phi_0, j_0)(H)$, and the claim follows. $\square$

## 2. Nonnegative curvature on homogeneous bundles

The purpose of this section is to prove Theorems E and G of the introduction.

As in [Ch1] we will construct nonnegative curvature metrics on $M = D(B_-) \cup_E D(B_+)$ (cf. 1.2) with the additional property that the common boundary $E = \partial D(B_-) = \partial D(B_+)$ is totally geodesic in $M$. This is a very strong restriction, which, by the soul theorem [CG], implies that also $B_\pm$ are to-tally geodesic. With this in mind, all we have to do is to construct $G$-invariant nonnegative curvature metrics on the bundles $D(B_\pm) = G \times_{K^\pm} D^{\ell_\pm+1}$ (cf. (1.1)), that agree on the common boundary $E = G/H = G \times_{K^\pm} S^{\ell_\pm}$ and are product metrics near the boundary.

From the Gray-O'Neill curvature submersion formula (cf. [ON] or [Gr]), we know that the product metric of a left invariant, $\mathrm{Ad}\,(K)$-invariant metric of nonnegative curvature on $G$ with a $K$-invariant nonnegative curvature metric on $D^{\ell+1}$ (which is product near $S^\ell = \partial D^{\ell+1}$) induces a $G$-invariant nonnegative curvature metric on the quotient $G \times_K D^{\ell+1}$ (which is product near $G/H = G \times_K S^\ell = \partial(G \times_K D^{\ell+1})$). The difficulty in the above strategy is therefore, that in general the restriction of such metrics on $G \times_{K_-} D^{\ell_-+1}$ and on $G \times_{K_+} D^{\ell_++1}$ to $G/H = G \times_{K_-} S^{\ell_-} = G \times_{K_+} S^{\ell_+}$ are different.

Consider any closed Lie subgroups $H \subset K \subset G$ of a compact Lie group $G$, with Lie algebras $\mathfrak{h} \subset \mathfrak{k} \subset \mathfrak{g}$. Fix any left invariant, $\mathrm{Ad}\,(K)$-invariant Riemannian metric, $\langle\ ,\ \rangle$ on $G$ and let $\mathfrak{m} = \mathfrak{k}^\perp$ and $\mathfrak{p} = \mathfrak{h}^\perp \cap \mathfrak{k}$ relative to this metric. On $G/H$ and $K/H$ we get induced (submersed) $G$-, respectively $K$-invariant metrics which are also denoted by $\langle\ ,\ \rangle$. As usual we make the identifications $\mathfrak{p}+\mathfrak{m} \simeq T_H G/H$ and $\mathfrak{p} \simeq T_H K/H$ via action fields; i.e., $X+A \to (X+A)^*_H$ and $X \to X^*_H$ respectively.

The homogeneous space $G/H$ can be identified with the orbit space $G \times_K K/H$ of $G \times K/H$ by the $K$-action $(k, (g, \bar{k}H)) \to (gk^{-1}, k\bar{k}H)$. The identification is given by $gH \to K(g, H)$ with inverse $K(g, kH) \to gkH$. By $\sqrt{\lambda}K/H$ we mean $K/H$ endowed with the metric $\lambda\langle\ ,\ \rangle$, where $\lambda > 0$. In this terminology we have:

LEMMA 2.1. *The $G$-invariant metric $\langle\ ,\ \rangle_\lambda$ on $G/H$ induced from the product metric on $G \times \sqrt{\lambda}K/H$ via $G \times_K K/H \simeq G/H$ is determined by*

$$\langle\ ,\ \rangle_{\lambda_{|\mathfrak{m}}} = \langle\ ,\ \rangle_{|\mathfrak{m}} \quad and \quad \langle\ ,\ \rangle_{\lambda_{|\mathfrak{p}}} = \tfrac{\lambda}{\lambda+1}\langle\ ,\ \rangle_{|\mathfrak{p}}.$$

*Proof.* The vertical space (= tangent space to $K$-orbit) at $(1, H) \in G \times K/H$ is given by

$$T^v_{(1,H)} = \mathfrak{h} \times \{0\} + \{(-X, X^*_H) \mid X \in \mathfrak{p}\}.$$

Thus $(U, Y^*_H) \in T_{(1,H)}G \times K/H$, $U \in \mathfrak{g}, Y \in \mathfrak{p}$ is horizontal if and only if $U = Z + A \in \mathfrak{p} + \mathfrak{m} = \mathfrak{h}^\perp$ satisfies $-\langle Z, X\rangle + \lambda\langle Y, X\rangle = 0$ for all $X \in \mathfrak{p}$; i.e.

$$T^h_{(1,H)} = \mathfrak{m} \times \{0\} + \{(\lambda Y, Y^*_H) \mid Y \in \mathfrak{p}\}.$$

Now $(A, 0)$ projects to $A \in \mathfrak{m} \subset T_H G/H$ and $(\lambda Y, Y^*_H)$ projects to $(\lambda+1)Y \in \mathfrak{p} \subset T_H G/H$. In particular, the horizontal lift of $A \in T_H G/H$ to $(1, H)$ is $(A, 0)$, and $Y \in \mathfrak{p} \subset T_H G/H$ lifts to $\frac{1}{\lambda+1}(\lambda Y, Y^*_H)$. This proves the claim since the norms of these vectors are given by $\|(A, 0)\|^2 = \|A\|^2$ and $\|\frac{1}{\lambda+1}(\lambda Y, Y^*_H)\|^2 = (\frac{1}{\lambda+1})^2(\lambda^2\|Y\|^2 + \lambda\|Y\|^2) = \frac{\lambda}{\lambda+1}\|Y\|^2$. $\qquad\square$

As an immediate consequence of this lemma, we see that if $Q$ is a fixed bi-invariant metric on $G$, and we choose $\langle \, , \, \rangle$ above as

$$(2.2) \qquad \langle \, , \, \rangle_{|\mathfrak{m}} = Q_{|\mathfrak{m}} \quad \text{and} \quad \langle \, , \, \rangle_{|\mathfrak{k}} = \tfrac{\lambda+1}{\lambda} Q_{|\mathfrak{k}}$$

then the metric on $G/H$ induced via $G \times_K \sqrt{\lambda} K/H$ as above, is the same as the one induced directly via $Q$. This is essentially the method that Cheeger used in [Ch1] to construct a nonnegatively curved metric on the connected sum of two projective spaces. The problem now, however, is that in general a metric like (2.2) on $G$ has some negative sectional curvature, as we will see, since $a = \frac{\lambda+1}{\lambda} > 1$.

We need to work in a slightly more general context. As before $G$ is a compact Lie group and $\mathfrak{k} \subset \mathfrak{g}$ a subalgebra. Let $K \subset G$ be the (immersed) Lie subgroup of $G$ with Lie algebra $\mathfrak{k}$; i.e. $K$ need not be compact. As before let $Q$ be a fixed bi-invariant metric on $\mathfrak{g}$ and $a > 0$. Define

$$(2.3) \qquad Q_{a|\mathfrak{m}} = Q_{|\mathfrak{m}} \quad \text{and} \quad Q_{a|\mathfrak{k}} = aQ_{|\mathfrak{k}}$$

and denote again by $Q_a$ also the corresponding left and $\mathrm{Ad}\,(K)$ invariant metric on $G$. We need the following curvature formulas for this left invariant metric (see, e.g., [Es] and [DZ] for special cases).

PROPOSITION 2.4. *For any $a > 0$ let $R^a$ be the curvature tensor of the metric $Q_a$ defined in (2.3). Then for any $A, B \in \mathfrak{m}$ and $X, Y \in \mathfrak{k}$ we have*

$$Q_a(R^a(A + X, B + Y)(B + Y), A + X)$$
$$= \tfrac{1}{4} \|[A, B]_{\mathfrak{m}} + a[X, B] + a[A, Y]\|_Q^2 + \tfrac{1}{4} \left\| [A, B]_{\mathfrak{k}} + a^2[X, Y] \right\|_Q^2$$
$$+ \tfrac{1}{4} a(1-a)^3 \|[X, Y]\|_Q^2 + \tfrac{3}{4}(1-a) \|[A, B]_{\mathfrak{k}} + a[X, Y]\|_Q^2$$

*where subscripts denote components. In particular, $(G, Q_a)$ has nonnegative curvature whenever $0 < a \le 1$, or if $\mathfrak{k}$ is abelian and $a \le \frac{4}{3}$.*

*Proof.* For $a = 1$ this is the well-known formula for the sectional curvature of a bi-invariant metric. For $a \ne 1$, we claim that $Q_a$ is a submersed metric. Indeed, on $G \times K$ consider the bi-invariant (semi-) Riemannian metric induced from $\langle \, , \, \rangle = Q \times bQ_{|_k}$ ($b$ negative allowed) on $\mathfrak{g} \times \mathfrak{k}$. When $b = \frac{a}{1-a}$ we claim that the map $G \times K \to G$, $(g, k) \mapsto gk$ is a (semi-) Riemannian submersion. In fact this can be viewed as a special case of (2.1) above, when $H$ is trivial, by noticing that in this case the vertical space given by

$$T_{(1,1)}^v = \{(-X, X) \mid X \in \mathfrak{k}\} \subset T_{(1,1)} G \times K$$

is nondegenerate since $b \ne -1$. (This would not be true in the general case where $\mathfrak{h} \ne \{0\}$.) The rest of the argument in (2.1) carries over verbatim and we see that the submersed metric on $G$ is scaled by $\frac{b}{b+1} = a$ in the $\mathfrak{k}$-direction.

To compute the sectional curvature, we use the Gray-O'Neill formula. Consider a 2-plane in $T_1G = \mathfrak{g}$ spanned by $A + X$ and $B + Y$ as in (2.4). The corresponding horizontal lifts to $T_{(1,1)}G \times K$ are $(A + \frac{b}{b+1}X, \frac{1}{b+1}X)$ and $(B + \frac{b}{b+1}Y, \frac{1}{b+1}Y)$, respectively. Moreover, when extending the $G$-coordinate to left-invariant vector fields and the $K$-coordinate to right-invariant vector fields, the resulting fields are easily seen to be horizontal. The Gray-O'Neill formula then yields:

$$Q_a(R^a(A + X, B + Y)(B + Y), A + X) = \alpha + \tfrac{3}{4}\beta,$$

where

$$\begin{aligned}
\alpha &= \langle R^{G \times K}\left((A + \tfrac{b}{b+1}X, \tfrac{1}{b+1}X), (B + \tfrac{b}{b+1}Y, \tfrac{1}{b+1}Y)\right) \\
&\qquad (B + \tfrac{b}{b+1}Y, \tfrac{1}{b+1}Y), (A + \tfrac{b}{b+1}X, \tfrac{1}{b+1}X)\rangle_{\mathfrak{g} \times \mathfrak{k}}
\end{aligned}$$

and

$$\beta = \left\|\left[(A + \tfrac{b}{b+1}X, \tfrac{1}{b+1}X), (B + \tfrac{b}{b+1}Y, \tfrac{1}{b+1}Y)\right]^v\right\|^2_{\mathfrak{g} \times \mathfrak{k}}.$$

Now

$$\begin{aligned}
\alpha &= \tfrac{1}{4}\left\|[A + \tfrac{b}{b+1}X, B + \tfrac{b}{b+1}Y]\right\|^2_Q + \tfrac{1}{4}b\left\|[\tfrac{1}{b+1}X, \tfrac{1}{b+1}Y]\right\|^2_Q \\
&= \tfrac{1}{4}\left\|[A,B]_{\mathfrak{m}} + \tfrac{b}{b+1}([X,B] + [A,Y])\right\|^2_Q \\
&\qquad + \tfrac{1}{4}\left\|[A,B]_{\mathfrak{k}} + (\tfrac{b}{b+1})^2[X,Y]\right\|^2_Q + \tfrac{1}{4}b(\tfrac{1}{b+1})^4\left\|[X,Y]\right\|^2_Q
\end{aligned}$$

where we have used $[\mathfrak{m}, \mathfrak{k}] \subset \mathfrak{m}$ and $[\mathfrak{k}, \mathfrak{k}] \subset \mathfrak{k}$. In terms of $a = \frac{b}{b+1}$ (and hence $1 - a = \frac{1}{b+1}$ and $b = \frac{a}{1-a}$) we have

$$\begin{aligned}
\alpha &= \tfrac{1}{4}\|[A,B]_{\mathfrak{m}} + a[X,B] + a[A,Y]\|^2_Q + \tfrac{1}{4}\left\|[A,B]_{\mathfrak{k}} + a^2[X,Y]\right\|^2_Q \\
&\qquad + \tfrac{1}{4}a(1-a)^3\|[X,Y]\|^2_Q.
\end{aligned}$$

Using the fact that for right-invariant vector fields $X^*, Y^*, [X^*, Y^*] = -[^*X, ^*Y] = -[X, Y]$ in terms of left-invariant vector fields, we get:

$$\begin{aligned}
\beta &= \left\|\left([A,B]_{\mathfrak{m}} + \tfrac{b}{b+1}[A,Y] + \tfrac{b}{b+1}[X,B] + [A,B]_{\mathfrak{k}} \right.\right. \\
&\qquad \left.\left. + (\tfrac{b}{b+1})^2[X,Y], -(\tfrac{1}{b+1})^2[X,Y]\right)^v\right\|^2_{\mathfrak{g} \times \mathfrak{k}} \\
&= \left\|\left([A,B]_{\mathfrak{k}} + (\tfrac{b}{b+1})^2[X,Y], -(\tfrac{1}{b+1})^2[X,Y]\right)^v\right\|^2_{\mathfrak{g} \times \mathfrak{k}},
\end{aligned}$$

since $\mathfrak{m} \times 0 \subset T^h$.

For $U, V \in \mathfrak{k}$ we have

$$(U, V) = \tfrac{1}{b+1}(b(U + V), U + V) + \tfrac{1}{b+1}(-(-U + bV), -U + bV)$$

and hence $(U, V)^v = \tfrac{1}{b+1}(U - bV, -U + bV)$. Moreover,

$$\|(U, V)^v\|^2_{\mathfrak{g} \times \mathfrak{k}} = (\tfrac{1}{b+1})^2 \left\{ \|U - bV\|^2_Q + b\| - U + bV\|^2_Q \right\} = \tfrac{1}{b+1}\|U - bV\|^2_Q.$$

This yields

$$\begin{aligned}
\beta &= \tfrac{1}{b+1} \left\| [A,B]_{\mathfrak{k}} + (\tfrac{b}{b+1})^2 [X,Y] + b(\tfrac{1}{b+1})^2 [X,Y] \right\|^2_Q \\
&= (1 - a) \left\| [A,B]_{\mathfrak{k}} + a[X,Y] \right\|^2_Q
\end{aligned}$$

which completes the proof of the curvature formula.

If $a \leq 1$ all terms are nonnegative. For $a > 1$ the latter two terms will be negative in general. However, if $\mathfrak{k}$ is abelian the formula reduces to

$$\begin{aligned}
Q_a(R^a(A + X, B + Y)(B + Y), A + X) \\
= \tfrac{1}{4} \|[A,B]_{\mathfrak{m}} + a[X,B] + a[A,Y]\|^2_Q + (1 - \tfrac{3}{4}a) \|[A,B]_{\mathfrak{k}}\|^2_Q
\end{aligned}$$

which is nonnegative for $a \leq 4/3$ as claimed. $\qquad \square$

*Remark* 2.5.  In general there are two planes with strictly negative curvature on $(G, Q_a)$ for any $a > 1$ arbitrarily close to 1. Indeed, one can usually easily find two planes spanned by $A + X$ and $B + Y$ with $[A, B] = -a^2[X, Y]$ and $[X, B] + [A, Y] = 0$ which will have negative sectional curvature if $[X, Y] \neq 0$.

We are now ready to prove the main result of this section.

THEOREM 2.6.  *Suppose $G$ is a connected, compact Lie group and $H \subset K \subset G$ are closed subgroups with $K/H = S^1$. Then for any bi-invariant metric on $G$, there is a $G$-invariant nonnegatively curved metric on $G \times_K D^2$ which is a product near the boundary $G/H = G \times_K S^1$, and so that the metric restricted to $G/H$ is induced from the given bi-invariant metric on $G$.*

*Proof.* Fix a bi-invariant metric $Q$ on $\mathfrak{g}$ and let $\mathfrak{m} = \mathfrak{k}^\perp, \mathfrak{p} = \mathfrak{h}^\perp \cap \mathfrak{k}$ as before. By assumption, $\mathfrak{p}$ is 1-dimensional and is hence an abelian subalgebra of $\mathfrak{g}$. Moreover, if $\bar{H} \subset H$ is the ineffective kernel of the $K$-action on $S^1 = K/H = (K/\bar{H})/(H/\bar{H})$ we have $\bar{\mathfrak{h}} = \mathfrak{h}$ since the isotropy group of an effective action on $S^1$ is finite. Since $\bar{H}$ is normal in $K$, $\mathfrak{h}$ and hence $\mathfrak{p}$ is preserved by $\mathrm{Ad}\,(K)$. This implies that the metric $\bar{Q}_a$ on $G$ defined by

$$\bar{Q}_{a|\mathfrak{m}} = Q_{|\mathfrak{m}}, \; \bar{Q}_{a|\mathfrak{p}} = aQ_{|\mathfrak{p}} \quad \text{and} \quad \bar{Q}_{a|\mathfrak{h}} = Q_{|\mathfrak{h}}$$

is $\mathrm{Ad}\,(K)$-invariant. Since $\mathfrak{p}$ is a subalgebra, this metric can also be viewed as in (2.3) with $\mathfrak{k} = \mathfrak{p}$ and $K = \exp(\mathfrak{p})$ since in (2.3) we allowed $K$ to be a noncompact Lie group. Hence (2.4) implies that the metric $\bar{Q}_a$ on $G$ has nonnegative sectional curvature if $a \leq 4/3$.

Let $K/H = S^1$ be equipped with the metric induced from $\bar{Q}_{a_{|\mathfrak{k}}}$. By Lemma (2.1), the metric on $G \times_K \sqrt{\lambda}S^1 \simeq G/H$ is then given by $Q$ on $\mathfrak{m}$ and $\frac{\lambda}{\lambda+1}aQ$ on $\mathfrak{p}$. Now pick, e.g., $a = 4/3$, $\lambda = 3$ and a $K$-invariant metric on $D^2$ with nonnegative curvature, which is a product near the boundary circle $\partial D^2 = S^1$, and on $S^1$ is the metric $\sqrt{3}S^1 = \sqrt{3}K/H$ from above.

The quotient metric on $G \times_K D^2$ induced from the product metric on $G$ and on $D^2$ has all the desired properties claimed in (2.6) .                    $\square$

From (2.6) and the discussion in the beginning of this section, it follows immediately that we can construct nonnegatively curved metrics on each half $G \times_{K_\pm} D^2$, matching smoothly near $\partial(G \times_{K_-} D^2) \simeq G/H \simeq \partial(G \times_{K_+} D^2)$ to yield $G$-invariant metrics on $M = G \times_{K_-} D^2 \cup_E G \times_{K_+} D^2$ with nonnegative curvature. This finishes the proof of Theorem E.

*Remark* 2.7.   For a metric on $D^2$ we can choose a rotationally symmetric metric $dt^2 + f(t)^2 d\theta^2$, where $f$ is a concave function which is odd with $f'(0) = 1$ in order to guarantee smoothness of the metric. Suppose $K/H = (S^1, Q)$ is a circle of length $2\pi r$. Then the induced metric on the principal orbit $G/H$ at $c(t)$ (where $t = 0$ corresponds to the singular orbit $G/K$) can be described as $G \times_K \frac{f(t)}{r\sqrt{a}}K/H$ which, using (2.1), is then given by $Q$ on $\mathfrak{m}$ and $\frac{f^2 a}{f^2 + ar^2}Q$ on $\mathfrak{p}$. Hence we need to choose a $t_0$ such that $f^2(t) = \frac{ar^2}{a-1}$, for $t \geq t_0$. Notice that the larger the radius $r$ is, or if we choose $1 < a \leq 4/3$ close to 1, the larger $t_0$ needs to be, and hence the diameter of $M$ will be large.

*Remark* 2.8.   In the case where a nonregular orbit is exceptional, i.e., is a hypersurface, one can just choose the bi-invariant metric on $G$ itself to induce a metric on the disc bundle $G \times_K D^1$, which then has the same properties as in Theorem 2.6. Hence one obtains a nonnegatively curved metric on every cohomogeneity one manifold with nonregular orbits of codimension $\leq 2$.

We point out that there are many cohomogeneity one manifolds with nonnegative curvature, whose singular orbits have codimension bigger than 2. One large class is the linear cohomogeneity one actions on round spheres $S^n(1)$, classified in [HL], and characterized as the isotropy representations of compact rank two symmetric spaces. There are also many isometric cohomogeneity one actions on compact symmetric spaces with their natural metric of nonnegative curvature, recently classified in [Ko] in the irreducible case. In almost all of these examples, none of the principal orbits are totally geodesic. The difficulty in proving the conjecture that every cohomogeneity one manifold carries a metric with nonnegative curvature may lie in that one needs a better understanding of how to glue the two halves together without making the middle totally geodesic.

One particularly intriguing class of cohomogeneity one manifolds are the $2n - 1$-dimensional Brieskorn varieties defined by the equations

$$z_0^d + z_1^2 + \cdots z_n^2 = 0 \quad , \quad |z_0|^2 + \cdots |z_n|^2 = 1.$$

For $n$ odd and $d$ odd, they are homeomorphic to spheres, and if in addition $d \equiv \pm 1 \mod 8$, they are diffeomorphic to spheres, whereas for $d \equiv \pm 3 \mod 8$, they are diffeomorphic to the Kervaire sphere. The Kervaire sphere is an exotic sphere if $2n - 1 \equiv 1 \mod 8$ or more generally if $n + 1$ is not a power of 2. As discovered in [HH], the Brieskorn variety carries a cohomogeneity one action by $\mathrm{SO}(2)\mathrm{SO}(n)$ defined by $(e^{i\theta}, A)(z_0, \cdots, z_n) = (e^{2i\theta}z_0, e^{id\theta}A(z_1, \cdots, z_n)^t)$. This action was examined in detail in [BH], where they showed that the group picture (for $d$ odd) is given by $K_- = \mathrm{SO}(2) \times \mathrm{SO}(n-2), K_+ = \mathrm{O}(n-1), H = \mathbb{Z}_2 \times \mathrm{SO}(n-2)$, with embeddings given by

$$(e^{i\theta}, A) \in K_- \subset \mathrm{SO}(2)\mathrm{SO}(2)\mathrm{SO}(n-2) \to (e^{i\theta}, R(d\theta), A)$$

with $R(d\theta)$ a rotation by angle $d\theta$, $A \in K_+ \to (\det(A), (\det(A), A))$, and $(\varepsilon, A) \in \mathbb{Z}_2 \times \mathrm{SO}(n-2) = H \to (\varepsilon, (\varepsilon, \varepsilon, A))$. In particular, one obtains a different action for each odd $d$, and the nonprincipal orbits have codimension 2 and $n - 1$.

In the special case $n = 3$ and $d$ odd, where these actions define a cohomogeneity one action on $S^5$, they were first discovered by E. Calabi, who also observed that they descend to cohomogeneity one actions on the homotopy projective spaces $S^5/Z_2$, where $Z_2$ is the element $-\mathrm{id} \in \mathrm{SO}(2)$. In [Lo] it was shown that this homotopy projective space contains four (oriented) diffeomorphism types, according to $d \equiv 1, 3, 5, 7 \mod 8$, and two homeomorphism types, according to $d \equiv \pm 1, \pm 3 \mod 8$. Notice that it is not known if all of the exotic $\mathbb{R}P^5$'s admit any orientation-reversing diffeomorphisms. Hence it is conceivable that two of the exotic differentiable structures are the same. In any case, each of the possible differentiable structures on $\mathbb{R}P^5$ carries infinitely many cohomogeneity one actions by $\mathrm{SO}(2)\mathrm{SO}(3)$, and since the codimension of the singular orbits in this case are both equal to two, they all admit an invariant metric with nonnegative sectional curvature by Theorem E. One easily shows that the effective group picture is given by $G = \mathrm{SO}(2)\mathrm{SO}(3)$, $K_- = \mathrm{SO}(2)$ with embedding $e^{i\theta} \to (e^{2i\theta}, (R(d\theta), \mathrm{id}))$, $K_+ = \mathrm{O}(2)$ with embedding $A \to (1, (\det(A), A))$ and $H = Z_2 = \langle (1, \mathrm{diag}\,(-1, -1, 1)) \rangle$.

To finish the proof of Theorem G, we need to show that these metrics are never isometric to each other. For this we first note that if the action of $\mathrm{SO}(2)\mathrm{SO}(3)$ extends to a transitive action, then it must be linear and hence corresponds to the case $d = 1$ which is the well-known tensor product action. If $d > 1$, we will argue that $\mathrm{SO}(2)\mathrm{SO}(3)$ is the identity component of the isometry group, and since the group actions are never conjugate to each other, the corresponding metrics cannot be isometric either. Notice that any isometries, besides the elements of SO(2)SO(3), must preserve the $G$ orbits and hence induce isometries of the homogeneous metrics on the principal orbits SO(2)(SO(3)/$Z_2$). One easily shows that for any invariant metric on this homogeneous space, any further isometries in the identity component come from right translations by $N^{\mathrm{SO(3)}}(Z_2)/Z_2$. But these right translations do not extend to $G/K_+$ and hence are not well defined on $M$. This finishes the proof of Theorem G.

## 3. Topology of principal bundles

In this section we discuss the proof of Theorem F from the introduction. First note that over $S^4$, every principal SO(2) bundle is trivial and well known obstruction theory implies that every $k$-dimensional vector bundle with $k > 4$ is the direct sum of a 4-dimensional bundle and a trivial bundle. Hence we only need to examine principal SO(3) and SO(4) bundles. This is also why Theorems E and F, together with the Gray-O'Neill submersion formula, imply Theorem B.

To employ the methods of Section 1 we begin by describing the well-known cohomogeneity one action by SO(3) on $S^4$ in a language that will be needed for our construction of principal bundles. Let

$$V = \{A \mid A \text{ a } 3 \times 3 \text{ real matrix with } A = A^t, \operatorname{tr}(A) = 0\}.$$

Then $V$ is a 5-dimensional vector space with inner product $\langle A, B \rangle = \operatorname{tr} AB$. SO(3) acts on $V$ via conjugation $g \cdot A = gAg^{-1}$ and this action preserves the inner product and hence acts on $S^4(1) \subset V$. Every point in $S^4(1)$ is conjugate to a matrix in $F = \{\operatorname{diag}(\lambda_1, \lambda_2, \lambda_3) \mid \sum \lambda_i = 0, \sum \lambda_i^2 = 1\}$ and hence the quotient space is 1-dimensional. The singular orbits $B_\pm$ consist of those matrices $A$ with two eigenvalues $\lambda_i$ the same, negative for $B_-$ and positive for $B_+$. Clearly, $F$ is a great circle in $S^4(1)$ that is orthogonal to all orbits and we can choose $x_- = \operatorname{diag}(2/\sqrt{6}, -1/\sqrt{6}, -1/\sqrt{6})$, $x_+ = \operatorname{diag}(1/\sqrt{6}, 1/\sqrt{6}, -2/\sqrt{6})$ and hence $K_- = \mathrm{S(O(1)O(2))}$, $K_+ = \mathrm{S(O(2)O(1))} \subset \mathrm{SO(3)}$. As long as $\lambda_1 > \lambda_2 > \lambda_3$ we obtain the principal isotropy group $H = \mathrm{S(O(1)O(1)O(1))} = \mathbb{Z}_2 \times \mathbb{Z}_2$. Notice that $B_-$ and $B_+$ are both Veronese surfaces in $S^4(1)$ which are antipodal to each other at distance $\pi/3$.

Next, we lift these groups into $S^3$ under the two-fold cover $S^3 = \mathrm{Sp}(1) \to \mathrm{SO(3)}$ which sends $q \in \mathrm{Sp}(1)$ into a rotation in the 2-plane $\operatorname{Im}(q)^\perp \subset \operatorname{Im}(\mathbb{H})$ with angle $2\theta$, where $\theta$ is the angle between $q$ and 1 in $S^3(1)$. After renumbering the coordinates, the group $K_-$ lifts to $\mathrm{Pin}(2) = \{e^{i\theta}\} \cup \{je^{i\theta}\}$ which we abbreviate to $e^{i\theta} \cup je^{i\theta}$. Similarly, $K_+$ lifts to $\mathrm{Pin}(2) = e^{j\theta} \cup ie^{j\theta}$, and $H = S(O(1)O(1)O(1)) \subset \mathrm{SO(3)}$ lifts to the quaternion group $Q = \{\pm 1, \pm i, \pm j, \pm k\}$.

Thus the group diagram for $S^4$ is

$$(3.1)$$

$$
\begin{array}{c}
S^3 \\
\diagup \quad \diagdown \\
e^{i\theta} \cup je^{i\theta} \qquad\qquad e^{j\theta} \cup ie^{j\theta} \\
\diagdown \quad \diagup \\
Q
\end{array}
$$

We are now in a position to construct principal SO(3) bundles over $S^4$. Since the second Stiefel-Whitney class $w_2$ of the principal bundle SO(3) $\to P^* \to S^4$ is zero, there exists a two-fold cover $P$ of $P^*$ such that $P \to S^4$ is a principal $S^3$ bundle. We first construct a cohomogeneity one action by $G$ on $P$, with $S^3 \subset G$, which then induces a cohomogeneity one action on $P^*$ since $P^* = P/\sigma$, with $\sigma = -1$ central in $S^3$, as long as $\sigma$ is also central in $G$.

Principal bundles $S^3 \to P \to S^4$ are classified by an element in $\pi_3(S^3) = \mathbb{Z}$ and hence by an integer $k$. Equivalently, we can consider the classifying map of the bundle $f \colon S^4 \to B_{S^3} = \mathbb{H}P^\infty$ and then $k = f^*(x)[S^4]$ where $x \in H^4(\mathbb{H}P^\infty, \mathbb{Z}) = \mathbb{Z}$ is the generator corresponding to $\mathbb{H}P^1 \subset \mathbb{H}P^\infty$. Hence we can also consider $k$ as the Euler class of the principal $S^3$ bundle, regarded as a sphere bundle over $S^4$, and evaluated on the fundamental class. Indeed the latter follows from the fact that the universal principal $S^3$ bundle over $\mathbb{H}P^\infty$ is the Hopf bundle with Euler class $x$. Throughout the rest of the paper we denote by $P_k \to S^4$ the principal $S^3$ bundle with Euler class $k$.

We can now use the $S^3$ cohomogeneity one action on $S^4$ in (3.1) and the main construction in (1.6) to arrive at the following group diagram:

$$(3.2)$$

$$
\begin{array}{c}
S^3 \times S^3 \\
\diagup \qquad\qquad \diagdown \\
(e^{ip_-\theta}, e^{i\theta}) \cup (j,j)(e^{ip_-\theta}, e^{i\theta}) \qquad\qquad (e^{jp_+\theta}, e^{j\theta}) \cup (i,i)(e^{jp_+\theta}, e^{j\theta}) \\
\diagdown \qquad\qquad \diagup \\
\triangle Q
\end{array}
$$

where $\triangle Q = \{\pm(1,1), \pm(i,i), \pm(j,j), \pm(k,k)\}$. In order for $H$ to be a subgroup of $K_\pm$, we need that $p_\pm \equiv 1 \bmod 4$ and then we get $K_\pm/H = S^1$. Hence (3.2) defines a cohomogeneity one manifold $P_{p_-,p_+}$. Notice that the action of $S^3 \times S^3$ is again ineffective, the effective version being $S^3 \times S^3/\pm(1,1) = \mathrm{SO}(4)$. As in (1.6), it now follows that $S^3 = S^3 \times 1$ acts freely on $P_{p_-,p_+}$ and that $P/S^3$ is a cohomogeneity one manifold as in (3.1) and hence equivariantly diffeomorphic to $S^4$. Thus we obtain a principal bundle

$$S^3 \to P_{p_-,p_+} \to S^4.$$

Since $\sigma = (-1, 1)$ is central in $S^3 \times S^3$, we also obtain a cohomogeneity one action by $SO(3) \times SO(3)$ on the principal $SO(3)$ bundle $P^* = P/(-1, 1) \to S^4$.

To identify the principal bundle, we prove:

PROPOSITION 3.3. *The principal $S^3$ bundle $P_{p_-,p_+} \to S^4$ is classified by* $k = (p_-^2 - p_+^2)/8$.

*Proof.* The Gysin sequence of the sphere bundle $S^3 \to P_k \to S^4$ yields that the nonzero cohomology groups of $P_k$ are: $H^0 = H^7 = \mathbb{Z}$ and $H^4(P_k, \mathbb{Z}) = \mathbb{Z}/|k|\mathbb{Z}$ if $k \neq 0$ and $H^3 = H^4 = \mathbb{Z}$ if $k = 0$. Hence we can recognize $|k|$ by computing the cohomology groups of $P_{p_-,p_+}$.

To do this in general for a cohomogeneity one manifold $M: H \subset \{K_-, K_+\} \subset G$, we use the Meyer-Vietoris sequence, where $U_\pm = D(B_\pm) = G \times_{K_\pm} D^{\ell_\pm + 1}$ deformation retracts to $B_\pm = G/K_\pm$ and $U_- \cap U_+ = G/H$. Hence we get a long exact sequence

$$(3.4) \quad \to H^{i-1}(B_-) \oplus H^{i-1}(B_+) \xrightarrow{\pi_-^* - \pi_+^*} \begin{aligned} &H^{i-1}(G/H) \to H^i(M) \\ &\to H^i(B_-) \oplus H^i(B_+) \to \end{aligned}$$

where $\pi_\pm$ are the projections of the sphere bundles $G/H = G \times_{K_\pm} S^{\ell_\pm} = \partial D(B_\pm) \to B_\pm = G/K_\pm$. Notice that in our case of (3.2) above, the restriction of the principal $S^3$ bundle $P_{p_-,p_+} \to S^4$ to the $S^3$ orbits $S^3/e^{i\theta} \cup je^{i\theta} \simeq \mathbb{R}P^2 \simeq S^3/e^{j\theta} \cup ie^{j\theta}$ and $S^3/Q$ in $S^4$ are all trivial, since the classifying space $\mathbb{H}P^\infty$ for principal $S^3$ bundles is 3-connected. Thus $B_\pm = G/K_\pm = S^3 \times \mathbb{R}P^2$ and $G/H = S^3 \times (S^3/Q)$ up to diffeomorphism. In particular we obtain: $H^3(B_\pm, \mathbb{Z}) = \mathbb{Z}$, $H^4(B_\pm, \mathbb{Z}) = 0$, and $H^3(G/H, \mathbb{Z}) = \mathbb{Z} + \mathbb{Z}$, and the Meyer-Vietoris sequence (3.4) for $P = P_{p_-,p_+}$ becomes:

$$\begin{aligned} H^3(P) \to H^3(B_-) \oplus H^3(B_+) &= \mathbb{Z} + \mathbb{Z} \xrightarrow{\pi_-^* - \pi_+^*} H^3(G/H) \\ &= \mathbb{Z} + \mathbb{Z} \to H^4(P) \to 0. \end{aligned}$$

In order to compute $H^4(P)$, we need to compute the cokernel of $\pi_-^* - \pi_+^*$. In our case, this cokernel is determined by the determinant of $\pi_-^* - \pi_+^*$. If the determinant is equal to 0, then $H^4(P) = \mathbb{Z}$, and if it is nonzero then $H^4(P)$ is a cyclic group with order the absolute value of the determinant. Consider the commutative diagram:

$$(3.5) \quad \begin{array}{ccc} S^3 \times S^3 & \xrightarrow{\tau_\pm} & S^3 \times S^3/K_\pm^o \\ \Big\downarrow{\eta} & & \Big\downarrow{\mu_\pm} \\ S^3 \times S^3/H & \xrightarrow{\pi_\pm} & S^3 \times S^3/K_\pm \end{array}$$

where $K_\pm^o$ are the identity components of $K_\pm$. The maps $\mu_\pm$ are two-fold covers and as before it follows that $S^3 \times S^3/K_\pm^o = S^3 \times (S^3/e^{i\theta}) = S^3 \times S^2$, and since $S^3 \times S^3/K_\pm = S^3 \times RP^2$, it follows that $\mu_\pm^*\colon H^3(G/K_\pm) \to H^3(G/K_\pm^o)$ is an isomorphism. The map $\eta$ is an 8-fold cover and if we write $S^3 \times S^3/H = S^3 \times (S^3/Q)$, then $\eta^*$ in $H^3$ is an isomorphism on the first factor and multiplication by 8 on the second. Hence $\eta^*\colon H^3(S^3 \times S^3/H, \mathbb{Z}) = \mathbb{Z}^2 \to H^3(S^3 \times S^3, \mathbb{Z}) = \mathbb{Z}^2$ has determinant 8 and therefore $\det(\pi_-^* - \pi_+^*) = \det(\tau_-^* - \tau_+^*)/8$. It remains to determine the induced map in $H^3$ for the $S^1$ bundle $\tau_\pm$. For this purpose we consider the following commutative diagram of fibrations, where we drop the $\pm$ index for the moment (see e.g. [WZ, p. 228]).

$$(3.6) \qquad \begin{array}{ccccccc}
S^1 & \longrightarrow & S^3 \times S^3 & \xrightarrow{\tau_\pm = \tau} & S^3 \times S^3/K^o & \xrightarrow{\rho_1} & B_{S^1} \\
\downarrow & & \downarrow{\scriptstyle \mathrm{id}} & & \downarrow & & \downarrow{\scriptstyle r} \\
S^1 \times S^1 & \longrightarrow & S^3 \times S^3 & \xrightarrow{h} & S^2 \times S^2 & \xrightarrow{\rho_2} & B_{S^1} \times B_{S^1}
\end{array}$$

coming from the $S^1$ bundle $\tau$ and the $S^1 \times S^1$ bundle $h$ (product of Hopf bundles). If we let $H^*(B_{S^1}) = \mathbb{Z}[s]$ and $H^*(B_{S^1} \times B_{S^1}) = \mathbb{Z}[t_1, t_2]$, then $r^*(t_1) = ps, r^*(t_2) = s$ since the inclusion $S^1 \to S^1 \times S^1$ is given by $e^{i\theta} \to (e^{ip\theta}, e^{i\theta})$. If we set $H^*(S^3 \times S^3) = \Lambda(u, v)$, then the only nonzero differentials in the spectral sequence for $\rho_2$ are $d_2(u) = t_1^2, d_2(v) = t_2^2$. By naturality the differentials in the spectral sequence for $\rho_1$ are given by $d_2(u) = p^2 s^2, d_2(v) = s^2$ and hence a generator 1 in $H^3(S^3 \times S^3/K^o)$ goes to $(-u, p^2 v)$ under $\tau^*$. Thus $\tau_-(1) = (-u, p_-^2 v), \tau_+(1) = (-u, p_+^2 v)$ and the matrix of $\tau_-^* - \tau_+^*$ is given by:

$$\begin{pmatrix} -1 & 1 \\ p_-^2 & -p_+^2 \end{pmatrix}$$

which implies that $|k| = |p_-^2 - p_+^2|/8$.

Next we will show that $k = \pm(p_-^2 - p_+^2)/8$ with a fixed choice of sign, that is, the sign does not depend on $p_-, p_+$. For this, consider the manifolds $P_{p_-, p_+}^7$ and $P_{p_+, p_-}^7$. We claim that the Euler class of the corresponding $S^3$ bundles differ by a sign. First note that the antipodal map $-\mathrm{id}\colon S^4 \to S^4$ interchanges the two halves of $S^4$ relative to the decomposition (3.1). Since it is orientation-reversing the Euler class of the pullback bundle $(-\mathrm{id})^* P_{p_-, p_+}$ is the negative of $P_{p_-, p_+}$. Moreover, $(-\mathrm{id})^* P_{p_-, p_+}$ is a cohomogeneity one manifold with diagram as for $P_{p_-, p_+}$, except the roles of $i$ and $j$ are switched. Precomposing the $S^3 \times S^3$ action by $(A, A)\colon S^3 \times S^3 \to S^3 \times S^3$ where $A$ is the inner automorphism of $S^3$ given by $A(i) = j, A(j) = i$ and $A(k) = k^{-1}$ we see that $P_{p_+, p_-}$ and $(-\mathrm{id})^* P_{p_-, p_+}$ are equivariantly diffeomorphic. In particular, the Euler class of $P_{p_+, p_-}$ and $P_{p_-, p_+}$ have opposite signs.

To see which sign is the correct one (although this is not important for our main results), we need to compute the Euler class in one particular case.

For this one can take the well-known cohomogeneity one action by SO(4) on $S^7$ (see e.g. [TT] ) which is given by the representation $\frac{3}{\circ}\hat{\otimes}\frac{1}{\circ}$ (the isotropy representation of the rank 2 symmetric space $G_2/\text{SO}(4)$). This action preserves the Hopf fibration $S^3 \to S^7 \to S^4$ with Euler class 1. By computing the isotropy groups of this SO(4) action, one shows (cf. [GZ1]) that they are the same as the ones for $P_{-3,1}$ and hence $k = (p_-^2 - p_+^2)/8$ as claimed. $\qquad\square$

COROLLARY 3.7. *Every principal $S^3$, respectively SO(3), bundle over $S^4$ has a cohomogeneity one action by $G = \text{SO}(4)$, respectively $G = \text{SO}(3) \times \text{SO}(3)$, in fact in general several inequivalent ones.*

*Proof.* We only need to convince ourselves that every integer $k$ can be written as $(p_-^2 - p_+^2)/8$, where $p_\pm \equiv 1 \bmod 4$. Set $p_- = 4r + 1, p_+ = 4s + 1$ and hence $k = (r - s)(2r + 2s + 1)$. Then if we let $r = -s$, we get $k = -2s$, for $r = s + 1$ we get $k = 4s + 3$ and for $r = s - 1$ we get $k = -4s + 1$. These solutions can also be written in the following more convenient form: $(p_-, p_+) = (2k+1, -2k+1)$ if $k \equiv 0 \bmod 2$, $(p_-, p_+) = (-k-2, -k+2)$ if $k \equiv 1 \bmod 4$ and $(p_-, p_+) = (k+2, k-2)$ if $k \equiv 3 \bmod 4$. Hence every integer $k$ can be achieved, in general in several different ways. $\qquad\square$

*Remark* 3.8. For each value of $k \neq 0$ there exist only finitely many solutions of $k = (p_-^2 - p_+^2)/8 = (r - s)(2r + 2s + 1)$, which can all be described as follows: Set $m = r - s$ and $n = 2r + 2s + 1$. Then $k = nm$ with $n$ odd, and $r = (2m + n - 1)/4, s = (-2m + n - 1)/4$. Hence for each way of writing $k$ as a product $nm$ with $n$ odd (including sign changes for both $n$ and $m$), we get a solution for $p_-$ and $p_+$, if $r$ and $s$ are integers. Notice that if $k = 2^t$, then $m = 2^t, n = 1$; hence one only gets one solution: $p_- = 2k+1, p_+ = -2k+1$ and it is not hard to see that in all other cases, one obtains several solutions. Thus all principal $S^3$ bundles with $k \neq 2^t$ have several inequivalent cohomogeneity one actions by $G = \text{SO}(4)$. If, for example, $k = 105$, then the following is the complete set of eight solutions: $(p_-, p_+) = (29, 1), (-31, -11), (37, -23), (41, 29), (-47, 37), (73, -67), (-107, -103), (-211, 209)$.

If $k = 0$, i.e., on $P^7 = S^4 \times S^3$, we obtain infinitely many inequivalent cohomogeneity one actions corresponding to $p_- = p_+$.

Using the principal $S^3$ bundle $P_k$ with Euler class $k$, we can consider the associated 2-sphere bundle $M_k = P_k \times_{S^3} S^2 \to S^4$, where $S^3$ acts on $S^2$ via the two-fold cover $S^3 \to \text{SO}(3)$. This can also be described as $M_k = P_k/S^1$ with $S^1 \subset S^3$. We now observe the following interesting consequence of our results:

COROLLARY 3.9. *The total space of the $S^2$ bundles $M_k \to S^4$, which admit a metric with nonnegative sectional curvature, have the same integral cohomology groups as $\mathbb{C}P^3$, but distinct cohomology rings for $k \geq 2$.*

*Proof.* The Gysin sequence of the sphere bundle $S^2 \to M_k \to S^4$ yields that the nonzero cohomology groups $H^*(M_k, \mathbb{Z})$ are $H^0 = H^2 = H^4 = H^6 = \mathbb{Z}$. From the Gysin sequence of the circle bundle $S^1 \to P_k \to M_k$, and $H^4(P_k, \mathbb{Z}) = \mathbb{Z}_k$, we get that if $x, y$ are the generators in $H^2(M_k, \mathbb{Z})$ and $H^4(M_k, \mathbb{Z})$, then $x^2 = ky$. Hence $M_k$ all have the same cohomology groups as $\mathbb{C}P^3$, but distinct cohomology rings, as long as $k \geq 2$. Notice that $M_k$ and $M_{-k}$ are diffeomorphic, $M_{\pm 1}$ is diffeomorphic to $CP^3$, and $M_0$ is diffeomorphic to $S^2 \times S^4$. $\qquad\square$

Next, we consider the case of principal $S^3 \times S^3$ bundles $P$ over $S^4$ and the corresponding principal SO(4) bundles $P^* \to S^4$ with $P^* = P/(-1, -1)$. These bundles are classified by elements of $\pi_3(S^3 \times S^3) = \pi_3(\mathrm{SO}(4)) = \mathbb{Z} \oplus \mathbb{Z}$ and hence by pairs of integers $(k, l)$. For this identification, we use the convention in [Mi]: To an element $(k, l)$ we associate the element in $\pi_3(\mathrm{SO}(4))$ given by $q \in S^3 \to (u \to q^k u q^l) \in \mathrm{SO}(4)$. Under the two-fold cover $S^3 \times S^3 \to \mathrm{SO}(4)$ given by $(q_1, q_2) \to (u \to q_1 u q_2^{-1})$ this corresponds to the element $q \to (q^k, q^{-l})$ in $\pi_3(S^3 \times S^3)$. Another way to describe these integers is as follows: If we start with a principal $S^3 \times S^3$ bundle $P_{k,l} \to S^4$, then we obtain two principal $S^3$ bundles $P/S^3 \times 1$ and $P/1 \times S^3$ and these are now classified by their Euler class $-l$ and $k$.

To construct cohomogeneity one actions on these principal bundles, we start with the group diagram

$$(3.10)$$

$$\begin{array}{ccc} & S^3 \times S^3 \times S^3 & \\ & \diagup \qquad \diagdown & \\ (e^{ip_-\theta}, e^{iq_-\theta}, e^{i\theta}) \cup (j,j,j)K_-^0 & & (e^{jp_+\theta}, e^{jq_+\theta}, e^{j\theta}) \cup (i,i,i)K_+^0 \\ & \diagdown \qquad \diagup & \\ & \Delta Q & \end{array}$$

which defines a cohomogeneity one manifold $P_{p_-,q_-,p_+,q_+}^{10}$ as long as $p_\pm, q_\pm \equiv 1 \bmod 4$. $S^3 \times S^3 \times 1$ acts freely on it with quotient $S^4$; hence $P_{p_-,q_-,p_+,q_+}$ is a principal $S^3 \times S^3$ bundle over $S^4$. As such, it is classified by two integers $k$ and $l$ as above. The analogue of Proposition 3.3 is now

PROPOSITION 3.11. *The principal $S^3 \times S^3$ bundle $P_{p_-,q_-,p_+,q_+} \to S^4$ is classified by $k = (p_-^2 - p_+^2)/8$ and $l = -(q_-^2 - q_+^2)/8$. Hence every principal $S^3 \times S^3$, respectively SO(4) bundle over $S^4$ has a cohomogeneity one action by $G = S^3 \times S^3 \times S^3/\pm(1,1,1)$, respectively $G = \mathrm{SO}(4) \times \mathrm{SO}(3)$.*

*Proof.* The formula for $k$ and $l$ follows from (3.3) since the group diagram for $P/S^3 \times 1 \times 1$ is the $S^3 \times S^3$ cohomogeneity one picture for $P_{q_-,q_+}$ and hence $l = -(q_-^2 - q_+^2)/8$ and similarly $k = (p_-^2 - p_+^2)/8$.

As before it follows that for each $k, l$, there exist solutions $p_\pm, q_\pm$ to $k = (p_-^2 - p_+^2)/8$ and $l = -(q_-^2 - q_+^2)/8$ with $p_\pm, q_\pm \equiv 1 \mod 4$. If $k \ne 0, l \ne 0$, there are only finitely many solutions, and for $k \ne 0, l = 0$, i.e., on $P_k^7 \times S^3$, there exist infinitely many different cohomogeneity one actions.

The ineffective kernel of the $S^3 \times S^3 \times S^3$ action is $\pm(1, 1, 1)$; hence on $P^* = P/(-1, -1, 1)$ the effective action is by

$$S^3 \times S^3 \times S^3 / \langle (-1, -1, -1), (-1, -1, 1) \rangle = SO(4) \times SO(3). \qquad \square$$

We finally point out that among the linear cohomogeneity one actions on spheres [HL], only $S^2, S^3, S^4, S^5$ and $S^7$ admit cohomogeneity one actions where both singular orbits have codimension 2. Moreover in each case there is only one effective action, and the groups are $S^1, T^2, SO(3), SO(2)SO(3)$ and $SO(4)$ respectively. Among the nonlinear cohomogeneity one actions on spheres, there exist infinitely many such actions by $SO(2)SO(3)$ on $S^5$ (cf. Section 2 and [St2]).

We now explore the consequences to the existence of the $SO(4)$ action on $S^7$ for vector bundles and sphere bundles over $S^7$. Since $\pi_6(SO(k)) = 0$ for $k = 2, 5, 6, 7$ (see [Ja]) it follows that only principal $SO(3)$ and $SO(4)$ bundles over $S^7$ can be nontrivial, and both admit two-fold covers to principal $S^3$ and $S^3 \times S^3$ bundles. We first consider the case of principal $S^3$ bundles. As was mentioned in the proof of (3.3), the cohomogeneity one picture for the $SO(4)$ action on $S^7$ is the same as that for $P_{-3,1}$. Hence, if we apply the construction in Section 1 to obtain principal $S^3$ bundles over $S^7$, one is forced to consider the same cohomogeneity one picture as that for $P_{p_-, -3, p_+, 1}$. By Proposition 3.11, $P_{p_-, -3, p_+, 1} = P_{k,1}$ is a principal $S^3 \times S^3$ bundle over $S^4$, where $k = (p_-^2 - p_+^2)/8$. One can of course also argue directly, that for every principal $S^3 \times S^3$ bundle $P_{k,1}$ over $S^4$, we have $P_{k,1}/S^3 \times 1 = P_1 = S^7$ and hence $P_{k,1}$ can be regarded as a principal $S^3$ bundle over $S^7$. As such, it is classified by an element $r \in \pi_6(S^3) = \mathbb{Z}_{12}$ (see [Ja]), and it was shown in [Ri3] that $r = k(k+1)/2$ and hence each principal $S^3$ bundle over $S^7$ with $r = 0, 1, 3, 4, 6, 7, 9, 10$ can be written in the form $P_{k,1}$ in infinitely many ways. We thus obtain:

COROLLARY 3.12. *Eight of the 12 principal $S^3$ bundles over $S^7$, classified by $r = 0, 1, 3, 4, 6, 7, 9$ and $10$, admit infinitely many cohomogeneity one actions by $S^3 \times S^3 \times S^3 / \pm (1, 1, 1)$.*

As a consequence, the associated bundles over $S^7$ with fiber $S^2$ or $\mathbb{R}^3$ also carry infinitely many metrics with nonnegative curvature. Note, however, as was done in [Ri3], that the total space of the principal $S^3$ bundles over $S^7$ not achieved by (3.12), i.e., $r = 2, 5, 8, 11$ are diffeomorphic to the corresponding ones for $r = 10, 7, 4, 1 \equiv -2, -5, -8, -11 \mod 12$. In fact, they are simply the pull back of these bundles via the reflection $R$ in the equator $S^6$ in $S^7$. As a consequence the corresponding associated $S^2$ and $R^3$ bundles also have diffeo­morphic total spaces. Thus all 3-dimensional vectorbundles and corresponding sphere bundles over $S^7$ have complete metrics of nonnegative curvature.

Similarly, principal $S^3 \times S^3$ bundles over $S^7$ are classified by $(r, s) \in \mathbb{Z}_{12} \oplus \mathbb{Z}_{12}$, and it follows as in (3.11) that every such bundle, with $r, s = 0, 1, 3, 4, 6,$ $7, 9, 10$, admits infinitely many cohomogeneity one actions by $S^3 \times S^3 \times S^3 \times S^3$. As before the principal bundles with $r, s = 2, 5, 8, 11$ as well as the correspond­ing associated bundles with fiber $S^3$ or $R^4$ have total spaces diffeomorphic to the ones with $r, s = 0, 1, 3, 4, 6, 7, 9, 10$. Also notice that the eight bundles with $(r, s) = (0, a)$ or $(a, 0)$, $a = 2, 5, 8, 11$ clearly have nonnegative curvature since the principal bundles are products of $S^3$ with principal $S^3$ bundles over $S^7$. In summary:

COROLLARY 3.13. *All three-dimensional and* 88 *of the* 144 *four-dimen­sional vector bundles over* $S^7$, *as well as the corresponding sphere bundles, have metrics with nonnegative sectional curvature.*

We conclude this section by pointing out that previously known methods yield the following (see also [Ri4]):

PROPOSITION 3.14. *All vector bundles and all sphere bundles over* $S^5$ *admit complete metrics of nonnegative curvature.*

*Proof.* Since $\pi_4(\mathrm{SO}(3)) = \mathbb{Z}_2, \pi_4(\mathrm{SO}(4)) = \mathbb{Z}_2 \oplus \mathbb{Z}_2$, and $\pi_4(\mathrm{SO}(5)) = \mathbb{Z}_2$ there are only 1,3, respectively 1 nontrivial vector bundle among the 3-, 4-, respectively 5-dimensional vector bundles over $S^5$. We will show that the total space of each of the corresponding principal bundles is diffeomorphic to a Lie group, such that the principal action is by isometries in the bi-invariant metric. This implies that all vector bundles and sphere bundles over $S^5$ admit a metric with nonnegative curvature.

The tangent bundle of $S^5$ gives rise to the nontrivial element in $\pi_4(\mathrm{SO}(5))$ $= \mathbb{Z}_2$ and its principal bundle is $\mathrm{SO}(6) \to \mathrm{SO}(6)/\mathrm{SO}(5) = S^5$.

As before, we can replace the principal $\mathrm{SO}(3)$ and $\mathrm{SO}(4)$ bundles with prin­cipal $S^3$ and $S^3 \times S^3$ bundles respectively. The principal $S^3$ bundle $\mathrm{SU}(3) \to$ $\mathrm{SU}(3)/\mathrm{SU}(2) = S^5$ is the nontrivial element in $\pi_4(S^3) = \mathbb{Z}_2$.

The nontrivial bundle corresponding to $(1, 0) \in \mathbb{Z}_2 \times \mathbb{Z}_2 = \pi_4(S^3 \times S^3)$ is given by the action of $S^3 \times S^3 = \mathrm{SU}(2) \times \mathrm{SU}(2)$ on $\mathrm{SU}(3) \times \mathrm{SU}(2)$:

$$(\alpha, \beta)(A, B) = (A\alpha^{-1}, B\beta^{-1})$$

and similarly for $(0, 1)$.

The action $(\alpha, \beta)(A, B) = (A\alpha^{-1}, \alpha B\beta^{-1})$ represents the nontrivial ele­ment $(1, 1)$. Indeed, if we divide by the first $\mathrm{SU}(2)$ the map

$$\mathrm{SU}(3) \times \mathrm{SU}(2)/\mathrm{SU}(2) \to \mathrm{SU}(3) : (A, B) \to A$$

is an SU(2)- equivariant diffeomorphism of SU(2) principal bundles, and if we divide by the second SU(2), the map

$$\text{SU}(3) \times \text{SU}(2)/\text{SU}(2) \to \text{SU}(3) : (A, B) \to AB$$

is an SU(2) equivariant diffeomorphism, and hence both bundles are nontrivial.

□

## 4. Almost free SO(3) actions

As we have seen in Section 3, there are typically many different ways of representing the principal bundles discussed in this paper as cohomogeneity one manifolds. This will in general yield different induced actions on associated bundles, and will enable us, in particular, to prove Theorems C and D in the introduction.

Recall, that any $S^3$ bundle over $S^4$ is associated to a principal SO(4) bundle over $S^4$, which in turn is determined by its two-fold universal cover, a principal $S^3 \times S^3$ bundle over $S^4$. Each of these bundles are thus determined by a pair of integers $(k, l) \in \mathbb{Z} \times \mathbb{Z} = \pi_3(S^3 \times S^3) = \pi_3(\text{SO}(4))$, where we use the convention described in the previous section. For $(k, l) \in \mathbb{Z} \times \mathbb{Z}$ let $M_{k,l} \to S^4, P^*_{k,l} \to S^4, P_{k,l} \to S^4$ denote the corresponding $S^3$ bundle, principal SO(4) bundle and principal $S^3 \times S^3$ bundle respectively. In Section 3 we saw that for any choice of integers $p_\pm, q_\pm \equiv 1 \mod 4$, satisfying $k = (p^2_- - p^2_+)/8$ and $l = -(q^2_- - q^2_+)/8$ there is a cohomogeneity one action by $S^3 \times S^3 \times S^3$ on $P_{k,l}$ with diagram (3.10), which induces an effective action of SO(4) × SO(3) on $P^*_{k,l}$. The SO(4) subaction is the free principal action on $P^*_{k,l}$ and the subaction by SO(3) is a lift of the cohomogeneity one action of SO(3) on $S^4$. In particular SO(3) acts on the total space of every associated $S^3$ bundle taking fibers to fibers.

THEOREM 4.1.  *The* SO(3) *action on* $M_{k,l}$, *induced from* (3.10) *as described above, preserves the* $S^3$ *fibration* $M_{k,l} \to S^4$ *and has exactly the following orbit types*:

$$(1), (\mathbb{Z}_2), (D_2), (D_{\frac{|p_- + q_-|}{2}}), (D_{\frac{|p_- - q_-|}{2}}), (D_{\frac{|p_+ + q_+|}{2}}), \ and \ (D_{\frac{|p_+ - q_+|}{2}})$$

*where* $D_0$, *in this context, should be interpreted as both* SO(2) *and* O(2).

*Proof.* To compute the isotropy groups of this action, we apply the Isotropy Lemma 1.7 to the corresponding ineffective $S^3$ action on $M_{k,l} = P^*_{k,l} \times_{\text{SO}(4)} S^3 = P_{k,l} \times_{S^3 \times S^3} S^3$. In the latter description $S^3 \times S^3$ acts on $S^3$ via quaternion multiplication $(Q_1, Q_2) \cdot v = Q_1 v Q_2^{-1}$. The isotropy groups of this $S^3 \times S^3$ action on $S^3$ are $\triangle S^3 \subset S^3 \times S^3$ and conjugates thereof, i.e. the subgroups $S^3_a = \{(b, aba^{-1}) \mid b \in S^3\}$ for some $a \in S^3$.

We can now read off the isotropy groups of the $S^3$-action from (1.7). They are $\phi_\pm^{-1}(S_a^3)$ and $\phi_0^{-1}(S_a^3)$, where $\phi_0\colon Q \to S^3 \times S^3$ is the diagonal embedding and $\phi_\pm\colon \mathrm{Pin}\,(2) \to S^3 \times S^3$ are the homomorphisms determined by $\phi_-(e^{i\theta}) = (e^{ip_-\theta}, e^{iq_-\theta})$, $\phi_-(j) = (j,j)$ and $\phi_+(e^{j\theta}) = (e^{jp_+\theta}, e^{jq_+\theta})$, $\phi_+(i) = (i,i)$.

Clearly $\phi_0^{-1}(S_a^3) = \langle -1 \rangle = \mathbb{Z}_2$ unless $a \in \langle e^{i\theta} \rangle, \langle e^{j\theta} \rangle, \langle e^{k\theta} \rangle$, and in these cases $\phi_0^{-1}(S_a^3) = \langle i \rangle, \langle j \rangle, \langle k \rangle = \mathbb{Z}_4$, except when $a = \pm 1$, in which case $\phi_0^{-1}(S_a^3) = Q$.

Now consider those $e^{i\theta_1}, je^{i\theta_2} \in \mathrm{Pin}\,(2)$ such that $(e^{ip_-\theta_1}, e^{iq_-\theta_1})$ or $j(e^{ip_-\theta_2}, e^{iq_-\theta_2}) \in S_a^3$. If $a = e^{it}$, then $ae^{ip_-\theta_1}a^{-1} = e^{iq_-\theta_1}$ implies that $e^{i(p_--q_-)\theta_1} = 1$ and $aje^{ip_-\theta_2}a^{-1} = je^{iq_-\theta_2}$ implies that $e^{i(p_--q_-)\theta_2 - 2it} = 1$. Hence

$$\phi_-^{-1}(S_{e^{it}}^3) = \langle e^{2\pi i/(p_--q_-)}, je^{2ti/(p_--q_-)} \rangle \subset \mathrm{Pin}\,(2)$$

for $p_- \neq q_-$. In the case of $p_- = q_-$, we get

$$\phi_-^{-1}(S_{e^{it}}^3) = \{e^{i\theta}\} = S^1 \subset \mathrm{Pin}\,(2)$$

if $a = e^{it} \neq \pm 1$ and $\phi_-^{-1}(S_{\pm 1}^3) = \mathrm{Pin}\,(2)$.

If $a = je^{it}$, then $ae^{ip_-\theta_1}a^{-1} = e^{iq_-\theta_1}$ implies that $e^{i(p_-+q_-)\theta_1} = 1$ and $aje^{ip_-\theta_2}a^{-1} = je^{iq_-\theta_2}$ implies that $e^{i(p_-+q_-)\theta_2 - 2it} = 1$. Hence

$$\phi_-^{-1}(S_{je^{it}}^3) = \langle e^{2\pi i/(p_-+q_-)}, je^{2ti/(p_-+q_-)} \rangle$$

as the only possibility, since $p_- \neq -q_-$ when $p_-, q_- \equiv 1 \bmod 4$.

If $a \notin \{e^{it}\} \cup j\{e^{it}\}$, then the only $\theta_1$ with $ae^{ip_-\theta_1}a^{-1} = e^{iq_-\theta_1}$ is given by $\theta_1 = 0, \pi$; i.e., $e^{i\theta_1} = \pm 1$. Moreover, for generic $a$, the equation $aje^{ip_-\theta_2}a^{-1} = je^{iq_-\theta_2}$ has no solutions. For special values of $a$ (depending on $p_-$ and $q_-$) there are precisely two values of $\theta_2$ ($\theta_2$ and $\theta_2 + \pi$) with $aje^{ip_-\theta_2}a^{-1} = je^{iq_-\theta_2}$. Hence $\phi_-^{-1}(S_a^3) = \mathbb{Z}_2$ or $\mathbb{Z}_4$.

The groups $\phi_+^{-1}(S_a^3)$ are computed in exactly the same way. Finally, to obtain the isotropy groups of the effective action by $\mathrm{SO}(3)$, we only need to observe that under the two-fold cover $S^3 \to \mathrm{SO}(3)$ the images of $Z_4$, $Q$, $\langle e^{2\pi i/p}, je^{2ti/p} \rangle$ (for $p$ even), $\{e^{i\theta}\}, \mathrm{Pin}\,(2)$ are equal to $Z_2$, $D_2$, $D_{p/2}$, $\mathrm{SO}(2)$ and $\mathrm{O}(2)$ respectively. □

As pointed out in (3.7), for each $(k,l)$ with $k \neq 0, l \neq 0$, there are only finitely many solutions $(p_\pm, q_\pm)$ to the equations $k = (p_-^2 - p_+^2)/8$, $l = -(q_-^2 - q_+^2)/8$, when $p_\pm, q_\pm \equiv 1 \bmod 4$. As explained there also, one of these solutions can be written as $(p_-, p_+) = (2k+1, -2k+1), (-k-2, -k+2)$ or $(k+2, k-2)$ when $k \equiv 0 \bmod 2, k \equiv 1 \bmod 4$ or $k \equiv 3 \bmod 4$ respectively. Similarly $(q_+, q_-) = (2l+1, -2l+1), (-l-2, -l+2)$ or $(l+2, l-2)$ when $l \equiv 0 \bmod 2, l \equiv 1 \bmod 4$ or $l \equiv 3 \bmod 4$ respectively. If say $l = 0$, $(q_-, q_+) = (4n+1, 4n+1)$ is obviously a solution for all $n$.

We exhibit the isotropy groups, other than $(1), (\mathbb{Z}_2)$ and $D_2$, of these particular $\mathrm{SO}(3)$ actions on $M_{k,l}$ in the following table:

|          | $k$ even                                          | $k$ odd                                             |
|----------|---------------------------------------------------|-----------------------------------------------------|
| l even   | $D_{|k+l|}, \ D_{|k-l\pm1|}$                       | $D_{|k+2l\pm1|/2}, \ D_{|k-2l\pm3|/2}$               |
| l odd    | $D_{|2k+l\pm1|/2}, \ D_{|2k-l\pm3|/2}$             | $D_{|k-l\pm4|/2}, \ D_{|k+l|/2}$                     |
| l = 0    | $D_{|2n+1\pm k|}, \ D_{|2n\pm k|}$                 | $D_{|4n+3\pm k|/2}, \ D_{|4n-1\pm k|/2}$             |

Table 4.2. Isotropy groups

In particular, we get:

COROLLARY 4.3.    *Each of the manifolds $M_{k,0}$ admit infinitely many inequivalent almost free $\mathrm{SO}(3)$ actions preserving the fibration $M_{k,0} \to S^4$ and inducing the same cohomogeneity one action on $S^4$.*

*Remark* 4.4.    Notice that $M_{k,0}$ are also precisely those $M_{k,l}$ which can be regarded not only as $S^3$ bundles over $S^4$, but also as principal $S^3$ bundles. Indeed, the glueing map $q \to \{u \to q^k u\}$ for the bundle $M_{k,0}$ commutes with the right action by $S^3$ and hence $S^3$ acts freely on $M_{k,0}$ with quotient $S^4$. Thus $M_{k,0} = P_k$ and one can therefore also directly lift the $\mathrm{SO}(3)$ action on $S^4$ using the cohomogeneity one action on $P_k = P_{p_-,p_+}$ from Corollary 3.7. But this is an effective action of $S^3$ on $M_{k,0}$, instead of $\mathrm{SO}(3)$, and one easily sees that it also acts almost freely with isotropy groups the binary dihedral groups $\langle e^{2\pi i/p_\pm}, j \rangle$. These actions of $S^3$, finitely many for each $k$, commute with the free principal action of $S^3$, whereas the infinitely many almost free actions of $\mathrm{SO}(3)$ in Corollary 4.3 do not.

A more detailed version of Corollary 4.3 in the special case of the Hopf fibration $S^7 = M_{1,0} \to S^4$ is included in the following result, which also implies Theorem C in the introduction.

THEOREM 4.5.    *For each $n$ there is an action of $\mathrm{SO}(3)$ on $S^7$ which preserves the Hopf fibration $S^7 \to S^4$ and has exactly the following orbit types*:

$$(1), (\mathbb{Z}_2), (D_2), (D_{|2n-1|}), (D_{|2n|}), (D_{|2n+1|}), (D_{|2n+2|})$$

*where as before $(D_0)$ stands for $(\mathrm{SO}(2))$ and $(\mathrm{O}(2))$. In particular, for $n \neq 0, -1$ this action is almost free. Moreover the action does not extend to the disc $D^8$ if $n \neq 0, \pm1, \pm2$.*

*Proof.* The first part can just be read off from Theorem 4.1 and Table 4.2 when $k = 1, l = 0$. To prove that the actions do not extend to the disc, we use the work of Oliver in [Ol] concerning the structure of fixed point free $\mathrm{SO}(3)$ actions on discs. First, however, consider the case of an $\mathrm{SO}(3)$ action on $D^8$ with nonempty fixed point set $D^{\mathrm{SO}(3)} \neq \emptyset$. Since SO(3) has only irreducible representations in odd dimensions, it follows that the slice representation of SO(3) at a fixed point, restricted to SO(2) $\subset$ SO(3), has to have a fixed vector and hence $\dim D^{\mathrm{SO}(2)} > 0$. By Smith theory $D^{\mathrm{SO}(2)} \supset D^{\mathrm{SO}(3)}$ has the integral cohomology of a point (cf. e.g. [Br, Chap. III]). In particular any component of $D^{\mathrm{SO}(2)}$ with positive dimension has nonempty boundary, and $\partial D^{\mathrm{SO}(2)} = D^{\mathrm{SO}(2)} \cap \partial D = S^{\mathrm{SO}(2)}$. Thus, if an almost free action of SO(3) on $S^7$ extends to $D^8$, it cannot have fixed points.

Next, we will show, using [Ol], that any fixed point free action of SO(3) on $D^8$ has $(\mathbb{Z}_3)$ or $(D_3)$ among its orbit types on the boundary sphere $S^7$, which then proves our theorem. From Corollary 1 of [Ol] we know in particular that $D_3$ occurs as isotropy group for any SO(3) action on $D^8$ without fixed points. In fact, from Lemmas 1 and 3 in [Ol], it follows that the octahedral group $O \subset \mathrm{SO}(3)$ has an isolated fixed point in the interior of $D^8$ and that $D_3$ occurs as an isotropy group of the linear representation of $O$ at such a fixed point. In particular, $D_3$ is the isotropy of an interior point $p \in D^8$ and $\dim D^{D_3} > 0$. Again by Smith theory $D^{\mathbb{Z}_3} \supset D^{D_3}$ has the $\mathbb{Z}_3$-cohomology of a point. Thus each component of $D^{\mathbb{Z}_3}$ intersects $\partial D = S$ nontrivially in $S^{\mathbb{Z}_3}$.

Relative to an SO(3)-invariant metric on $D^8$, join $p \in D^{D_3} \subset D^{\mathbb{Z}_3}$ to a closest point $q \in \partial D^{\mathbb{Z}_3} = S^{\mathbb{Z}_3}$ inside $D^{\mathbb{Z}_3}$. In particular, $\mathbb{Z}_3 \subset \mathrm{SO}(3)_q$. But $\mathrm{SO}(3)_q$ also fixes the normal vector to the boundary at $q$, hence the minimal geodesic to $p$. Hence $\mathrm{SO}(3)_q \subset \mathrm{SO}(3)_p = D_3$, which implies that $\mathrm{SO}(3)_q = \mathbb{Z}_3$ or $D_3$. $\qquad\square$

*Remark* 4.6.  The group of linear symmetries of the Hopf fibration $S^7 \to S^4$ is given by $(\mathrm{Sp}(2) \times \mathrm{Sp}(1))/\mathbb{Z}_2$ where $\mathrm{Sp}(2)$ acts via matrix multiplication on $S^7 \subset \mathbb{H}^2$, and $\mathrm{Sp}(1)$ is the right Hopf action. The cohomogeneity one action of SO(3) on $S^4$ defines an embedding of a maximal SO(3) in SO(5), which under the two-fold cover $\mathrm{Sp}(2) \to \mathrm{SO}(5)$ lifts to a maximal $\mathrm{Sp}(1)_m \subset \mathrm{Sp}(2)$. Hence $(\mathrm{Sp}(1)_m \times \mathrm{Sp}(1))/\mathbb{Z}_2 = \mathrm{SO}(4) \subset (\mathrm{Sp}(2) \times \mathrm{Sp}(1))/\mathbb{Z}_2$, which also happens to be the cohomogeneity one action of SO(4) on $S^7$ with singular orbits of codimension two, projects to $\mathrm{SO}(3) \subset \mathrm{SO}(5)$ under the Hopf map. Thus there are two linear lifts of the cohomogeneity one action of SO(3) on $S^4$. One is the almost free action by $Sp(1)_m$ on $S^7$ with isotropy groups 1 and $\mathbb{Z}_3$, and the other is the action by $\mathrm{SO}(3) = \triangle\mathrm{Sp}(1)/\mathbb{Z}_2 \subset \mathrm{SO}(4)$, which has isotropy groups 1, $\mathbb{Z}_2$, $D_2$, $\mathrm{SO}(2)$, and $\mathrm{O}(2)$. In particular, none of the actions in (4.5) are linear, except $n = 0$, which is the latter one.

We now analyze the ramifications of (4.1) to the Milnor spheres. Recall that Milnor [Mi] showed that the Euler class of $M_{k,l} \to S^4$ is equal to $e = k + l$. Hence the Gysin sequence and Smale's solution of the Poincaré conjecture implies that $M_{k,l}$ is homeomorphic to $S^7$ if and only if $k + l = \pm 1$. By changing the orientation if necessary, we can assume that $k + l = 1$. For Theorems A and D, we also need the diffeomorphism classification of these homotopy spheres $M_{k,1-k}$ due to Eells and Kuiper [EK]. According to [EK], $M_{k,1-k}$ is oriented diffeomorphic to $M_{m,1-m}$ if and only if $k(k-1) \equiv m(m-1)$ mod 56, i.e. the oriented diffeomorphism class is given by $\frac{k(k-1)}{2}$ mod 28 in the group $\mathbb{Z}_{28}$ of exotic 7-spheres. Notice also that a change in sign in $\mathbb{Z}_{28}$ corresponds to a change in orientation of the manifold and hence the numbers 1 to 14 correspond to the 14 possible distinct diffeomorphism types of exotic spheres.

As was observed in [EK], $\frac{k(k-1)}{2}$ mod 28 takes on the following 16 values: 0,1,3,6,7,8,10,13,14,15,17,20,21,22,24,27 which via a change of orientation corresponds to the numbers 0,1,3,4,6,7,8,10,11,13,14 and hence 11 of the 15 different diffeomorphism types of topological 7-spheres fiber over $S^4$ with $S^3$ as fiber. Using Theorems E and F, this completes the proof of Theorem A in the introduction. In passing, we note that $M_{2,-1}$ generates the group $Z_{28}$ of all homotopy 7-spheres via connected sum.

It is now clear that the SO(3) actions considered in (4.1) (cf. Table 4.2) on $M_{k,1-k}$ and $M_{m,1-m}$ are, in general, different actions on the same homotopy sphere when $\frac{k(k-1)}{2} \equiv \frac{m(m-1)}{2}$ mod 28. To make this more concrete, we exhibit the following special cases. As pointed out in [EK, p. 102], $k(k-1) \equiv m(m-1)$ mod 56 if and only if $m \equiv k$ or $1-k$ mod 7 , and $m \equiv k$ or $1-k$ mod 8 . Choosing the special case $m \equiv k$ mod 56, we get:

COROLLARY 4.7. *Let $M_{k,1-k}$ be any of the homotopy 7-spheres considered above. Then for each integer $n$, $M_{k,1-k}$ supports an* SO(3) *action with the following orbit types*:

$$(1), (\mathbb{Z}_2), (D_2), (D_{|k+56n+1\pm1|/2}), (D_{|3(k+56n)-1\pm3|/2})$$

*if $k$ is even, and*

$$(1), (\mathbb{Z}_2), (D_2), (D_{|k+56n-2\pm1|/2}), (D_{|3(k+56n)-2\pm3|/2})$$

*if $k$ is odd.*

Of course, there are many more actions given by Theorem 4.1 on the exotic spheres $M_{k,1-k}$, most of which are almost free. Even for the standard sphere, we get many additional almost free actions, besides the ones described in Theorem C, whenever $M_{k,1-k}$ is diffeomorphic to $S^7$, i.e. for $k \equiv 0, 1$ mod 7 and $k \equiv 0, 1$ mod 8. They preserve a different fibration of $S^7$ by 3-spheres, but in this case, we get only finitely many actions for each fibration.

Not all of the actions in Table 4.2 and Corollary 4.7 are almost free. Indeed, in the case of $M_{k,l}$ with $k = -l = 2^r$ there exists no almost free lift preserving the fibration, since (3.8) implies that the only action obtained from
(4.1) is the one described in Table 4.2. For the homotopy spheres $M_{k,1-k}$,
the only actions in Table 4.2 which are not almost free occur in the case
of $k = 0, 1, -2, 3$. Of course for $k = 0, 1$, which corresponds to the Hopf
fibration, (4.5) gives rise to infinitely many almost free actions preserving the
fibration. For $k = -2, 3$, i.e. on $M_{-2,3} = M_{3,-2}$, (4.1) implies that there
exist one further action besides the one described in Table 4.2. It corresponds
to $(p_-, p_+) = (5, 1)$, $(q_-, q_+) = (-3, 5)$ and hence gives rise to an almost free
action with isotropy groups $1, Z_2, D_2, D_3, D_4$. Thus in the case of the homotopy
spheres $M_{k,1-k}$ there always exists at least one almost free action preserving
the fibration. This implies Theorem D in the introduction.

We also observe that the SO(3) actions on $M_{k,l}$ extend to an action of
O(3). For this just note that the element $-\mathrm{id} \in \mathrm{SO}(4)$ commutes with the
structure group and the SO(3) action on $P_{k,l}^*$ and hence induces an action on
the associated bundle.

All the actions in this section on $M_{k,l}$ are isometric actions with respect to
the nonnegatively curved metrics we constructed in Section 2. We now consider
the question whether these metrics can ever be isometric to each other, and
show that, at least in the case of the homotopy spheres $M_{k,1-k}$, this can almost
never be the case:

PROPOSITION 4.8.   *If $M_{k,1-k}$ and $M_{m,1-m}$ are diffeomorphic, then the
metrics of nonnegative curvature constructed on them in Section 2 can only
be isometric, if the corresponding isometric* SO(3) *actions are conjugate. In
particular, we obtain infinitely many such metrics on each $M_{k,1-k}$ which are
not isometric to each other.*

*Proof.* To see this, we use the result by E. Straume that the degree of
symmetry of any exotic 7-sphere is at most 4, see [St1, Theorem C]. In other
words the dimension of any compact Lie group $G$ that acts effectively on an
exotic 7-sphere is at most 4. Now fix a metric on $\Sigma = M_{k,1-k}$ such that one of
the above SO(3) actions is isometric and let $G \supset \mathrm{SO}(3)$ be the id-component
of its full isometry group. Following Straume, $G$ is either SO(3) or a finite
quotient of $\mathrm{SO}(3) \times \mathrm{SO}(2)$. In particular $G$ contains only one subgroup SO(3),
so if one of the other actions of SO(3) on $\Sigma$ is a subgroup of $G$, then the two
actions must be conjugate. Hence, if the SO(3) actions are inequivalent, the
corresponding metrics cannot be isometric and in fact have different isometry
groups.                                                                      □

We suspect that O(3) will always be the full isometry groups of the met-
rics we constructed on $M_{k,l}$; see the next sections for some comments on this
question.

## 5. Remarks and open problems

Recall the two steps in our approach to the Cheeger-Gromoll problem: (1) Any principal SO(k)-bundle over $S^4$ has a cohomogeneity one $G$-structure with SO($k$) $\subset G$ (and with singular orbits of codimension two). (2) Any cohomogeneity one $G$-manifold (with singular orbits of codimension two) admits a $G$-invariant metric with nonnegative curvature. As suggested in the introduction, it is plausible that any cohomogeneity one manifold supports an invariant metric with nonnegative curvature. This is just one of the reasons for the following challenging:

PROBLEM 5.1. *Which principal* SO($k$)-*bundles over* $S^n$ *with* $k \leq n$ *support a cohomogeneity one $G$-structure with* SO($k$) $\subset G$?

Note that for each of the ways of writing $S^n$ as a cohomogeneity one manifold (cf. [HL]), our construction in Section 1 will in general yield several candidates for such bundles. In some cases it will give rise to infinitely many such candidates, namely whenever $S^1$ is a normal subgroup of $K_-$ or $K_+$. One can further increase the flexibility of our construction by using subactions of the usual cohomogeneity one actions listed in [HL], which are still cohomogeneity one (see [St2] for a complete list), and by making the actions ineffective. Of particular interest here are of course SO(8)-bundles over $S^8$, since 4095 exotic 15-spheres can be presented as (linear) 7-sphere bundles over the 8-sphere (cf. [Sh] and [EK]).

In view of our examples, it would be interesting to study in more detail the topology of the principal $S^3 \times S^3$ bundles $P_{k,l} \to S^4$ and their associated 3-sphere bundles $M_{k,l} \to S^4$. As we observed before, in the case of the sphere bundles $M_{k,l}$, one can recover the Euler class $e = k + l$ from the torsion in $H^4$. But for the principal bundles $P_{k,l}$ the only nonzero cohomology groups $H^*(P_{k,l}, \mathbb{Z})$ are $H^0 = H^3 = H^7 = H^{10} = \mathbb{Z}, H^4 = \mathbb{Z}_{(k,l)}$. Indeed, in Section 3 we saw that among the total spaces $P_{k,1}$ there are at most seven diffeomorphism types, and, using [HiR], it follows that there are precisely seven diffeomorphism types. In the special case of the homotopy spheres $M_{k,1-k}$, the total space has been classified up to diffeomorphism in [EK]. It would be interesting to extend this classification:

PROBLEM 5.2. *Classify the manifolds* $P_{k,l}$ *and* $M_{k,l}$ *up to homotopy, homeomorphism and diffeomorphism type.*

See [JW] and [Tam] for a partial classification of $M_{k,l}$ up to homotopy and homeomorphism type. For general results about 2-connected 7-manifolds see [Wa] and [Wi]. Also notice that it was shown in [DW] that for the corresponding vector bundles $E_{k,l} \to S^4$, the total spaces are diffeomorphic if and only if they are isomorphic as vector bundles.

Since our manifolds $M_{k,l}$ have trivial $\pi_2$ and $\pi_3 = \mathbb{Z}_{k+l}$, one easily sees that if one of the manifolds $M_{k,l}$ is homotopy equivalent to a homogeneous space $G/H$, then $G/H = \mathrm{Sp}(2)/\mathrm{Sp}(1)$, where $\mathrm{Sp}(1)$ is one of the three possible embeddings of $\mathrm{Sp}(1)$ in $\mathrm{Sp}(2)$. Hence $G/H = \mathrm{Sp}(2)/\mathrm{Sp}(1) \times 1 = S^7$, $G/H = \mathrm{Sp}(2)/\triangle\mathrm{Sp}(1) = T_1 S^4$ or $G/H = \mathrm{Sp}(2)/\mathrm{Sp}(1) = \mathrm{SO}(5)/\mathrm{SO}(3)$ where $\mathrm{SO}(3) \subset \mathrm{SO}(5)$ is the maximal embedding given by the cohomogeneity one action of $\mathrm{SO}(3)$ on $S^4$. It was shown in [Be] that $B^7 = \mathrm{SO}(5)/\mathrm{SO}(3)$ carries a metric of positive sectional curvature and that $H^4(B^7, \mathbb{Z}) = \mathbb{Z}_{10}$. One easily shows that the first Pontryagin class of the tangent bundle of $B^7$ is equal to 6 times a generator in $H^4 = \mathbb{Z}_{10}$. Furthermore, in [Tam] it was shown that the first Pontryagin class of the tangent bundle of $M_{k,l}$ is equal to $\pm 4l$ times a generator in $H^4(M_{k,l}, \mathbb{Z}) = \mathbb{Z}_{k+l}$. Hence $B^7$ cannot be homeomorphic to a principal $S^3$ bundle over $S^4$; i.e. $l = 0$. But it would be interesting to know if it can be homeomorphic or diffeomorphic to a sphere bundle.

For the principal $S^3$ bundles $P_k$ over $S^4$, we have that $P_k$ is diffeomorphic to $P_{-k}$ and that $P_{\pm 1} = S^7$. Furthermore, as was observed in [Ri2], $P_{\pm 2} = T_1 S^4$ since on $T_1 S^4 = \mathrm{Sp}(2)/\triangle\mathrm{Sp}(1)$ one has the free action by $S^3$ given by left multiplication with $\mathrm{diag}\,(q, 1)$ and since $H^4(T_1 S^4, \mathbb{Z}) = \mathbb{Z}_2$, it follows that this principal $S^3$ bundle has $k = \pm 2$. Hence $P_{\pm 1}$ and $P_{\pm 2}$ are diffeomorphic to homogeneous spaces. From the above, it follows that all other $P_k$ are strongly inhomogeneous; i.e. they do not have the homotopy type of a homogeneous space, except possibly $P_{\pm 10}$ which is at least not homeomorphic to a homogeneous space.

For the associated 2-sphere bundles $M_k \to S^4$ considered in (3.9), it follows from [On, Theorem 6] that the only homogeneous spaces that have the same integral cohomology groups as $\mathbb{C}P^3$, are $\mathbb{C}P^3 = M_{\pm 1}$ itself and $S^2 \times S^4 = M_0$. Hence $M_k$ with $|k| \geq 2$ do not have the homotopy type of a homogeneous space.

By construction, the total space $P$, respectively $M$, of any principal bundle, respectively associated sphere or vector bundle considered in this paper, as well as the base $S^4$, is the union $D_- \cup D_+$ of two disc bundles with common boundary $S = \partial D_- = \partial D_+$. Moreover, relative to the metrics of nonnegative curvature on $D_- \cup D_+$, $S$ is totally geodesic. From the Cheeger-Gromoll soul-construction [CG] and Perelman's rigidity theorem [Pe] it follows in particular that there are two-planes with zero curvature at every point of $D_- \cup D_+$.

All cohomogeneity one $G$-manifolds considered in this paper have $G = S^3, S^3 \times S^3$ or $S^3 \times S^3 \times S^3$ (acting possibly ineffectively). For the metrics on $D_\pm = G \times_{K_\pm} D^2$, note that we can choose a fixed bi-invariant metric on $G$ scaled by a fixed constant (e.g. $a = 4/3$) in the $K_\pm$ direction, and on $D^2$ we choose a metric $dt^2 + f^2(t)d\theta^2$ with a fixed convex function $f$. Although $K_\pm$ and hence $Q_a$ depends on the particular example, it easily follows from (2.4) and the Gray-O'Neill submersion formula that there is a uniform bound $C$ for the curvatures of all principal bundles considered here; i.e. $0 \leq \sec(P) \leq C$. But, as explained in (2.7), there is no bound on the diameter, $\mathrm{diam}(P)$, since the length of the circles $K_{\pm}$ goes to infinity. It is also apparent that all examples have a uniform lower bound on their volumes; i.e., there is a $v > 0$ such that $\mathrm{vol}\,(D_- \cup D_+) \geq v$ since this is true of $\mathrm{vol}\,(S)$. Similar bounds for curvature, volume and diameter hold for the associated bundles and the base $S^4$.

All these examples complement Cheeger's classical finiteness theorem [Ch2] and recent finiteness theorems by Petrunin-Tuschmann [PT] and Tapp [Ta].

Since our examples of $S^3$ bundles over $S^4$ are 2-connected, they illustrate the sharpness of the following, even within the class of nonnegatively curved manifolds.

THEOREM 5.3 (Petrunin-Tuschmann).  *For each $n$ and $D$, $C > 0$ there exist only finitely many diffeomorphism types of simply connected compact Riemannian $n$-dimensional manifolds $M$, with $|\sec M| \leq C$, $\mathrm{diam}\, M \leq D$ and finite $\pi_2(M)$.*

In the special case where the lower curvature bound is a fixed positive number $\delta > 0$, the same conclusion was obtained simultaneously by Fang and Rong in [FR] (in that case the bound on $\mathrm{diam}\, M$ is automatic by the Bonnet-Myers theorem). Motivated by this, the following conjecture was proposed in [FR]:

CONJECTURE (Rong).  *For each $n$, there are at most finitely many 2-connected positively curved $n$-manifolds.*

Our examples show that this conjecture is false, if we replace positive with nonnegative curvature.

The following result was first obtained in [GW] in the special case where $\Sigma = S^n(1)$. It was then extended to arbitrary souls in [Ta].

PROPOSITION 5.4.  *For each $n$ and $C, D, V > 0$ and each metric on $\Sigma$ with $\mathrm{diam}\, \Sigma \leq D$ and $\mathrm{vol}\, \Sigma > V$, there exist only finitely many vector bundles $M$ over $\Sigma$ with a complete metric of nonnegative curvature, such that $\Sigma$ is the soul and such that the sectional curvatures of all 2-planes of $M$ either tangent to $\Sigma$ or normal to $\Sigma$ are bounded above by $C$.*

By the above remarks, in our examples of 3- and 4-dimensional vector bundles over $S^4$, we have bounds on the curvatures of $M$ and the volume of the soul $\Sigma$, which is always the zero section of the vector bundle and isometric to the metric on the base $S^4$. But the diameter of the base necessarily goes to infinity since the length of the circles $K_{\pm}/H = S^1$ goes to infinity.

We conclude our discussion with a few more remarks about the geometry and symmetry of our examples, in particular the principal $S^3 \times S^3$ bundles $P_{k,l} \to S^4$ and their associated sphere bundles $S^3 \to M_{k,l} \to S^4$ and vector bundles $\mathbb{R}^4 \to E_{k,l} \to S^4$. Much work has been done previously on trying to construct metrics with nonnegative or positive curvature on the total spaces $P_{k,l}, M_{k,l}$ and $E_{k,l}$. A natural approach is to consider Kaluza-Klein type metrics on the principal $L$ bundle $P$, where one chooses a principal connection to define the horizontal space, pulls back the metric from the base to the horizontal space, and defines the metric on the fiber to be a bi-invariant or left-invariant metric on $L$. This metric then also induces a metric on the associated sphere bundles and vector bundles. We call these metrics *connection-type metrics*. In all three cases, the fibers of the projection onto the base are totally geodesic and isometric to each other. The metrics of positive Ricci curvature on $P_{k,l}$ and $M_{k,l}$ constructed in [Na] and [Po] are exactly of this type. But for the construction of nonnegatively curved metrics this approach has been successful only in the case where the principal bundle is a Lie group or a homogeneous space. In [DR] it was shown that the only case in which the induced metric on $M_{k,l}$ has positive curvature, is when $k = 0, l = \pm 1$ or $k = \pm 1, l = 0$, i.e. when $M_{k,l} = S^7$. It would be interesting to know if nonnegatively curved connection-type metrics exist on the bundles $M_{k,l}$ with $kl \neq 0, 1$. The metrics in our examples are not of this type. In our case, we will show that the metrics on the $S^3$ fibers (as well as on the base $S^4$) are cohomogeneity one metrics, and that the fibers are not totally geodesic.

On $P_{k,l}$ we have the cohomogeneity one action by $G = S_1^3 \times S_2^3 \times S_3^3$ with the principal bundle action given by $S_1^3 \times S_2^3$ (but acting on the left on $P$). By construction, the projection $P_{k,l} \to S^4$ is a Riemannian submersion, with the horizontal distribution given by a principal connection, since the metric is $S_1^3 \times S_2^3$ invariant. The same follows for the associated bundles $M_{k,l} = P_{k,l} \times_{S_1^3 \times S_2^3} S^3$ and $E_{k,l} = P_{k,l} \times_{S_1^3 \times S_2^3} \mathbb{R}^4$. In all three cases, the metric on the base is given by the submersed metric on $P_{k,l}/S_1^3 \times S_2^3$, and as a cohomogeneity one metric on $S^4$ under the action of $S_3^3$, is described by three functions $f_1(t), f_2(t), f_3(t)$, the length of the three action fields $i^*, j^*, k^*$ along $c(t)$ with $i, j, k$ in the Lie algebra of $S_3^3$. Here $c(t)$ is a fixed geodesic perpendicular to all orbits, as in Section 1. Invariance of the metric under the isotropy action of $H = Q$ implies that these action fields must be orthogonal. It follows from Perelman's rigidity theorem that two of these functions are equal to 1, $f_2(t) = f_3(t) = 1$ on $D_-$ and $f_1(t) = f_3(t) = 1$ on $D_+$.

The metric on the fiber of $P_{k,l} \to S^4$ over the point $c(t)$ in $S^4$ is given by a left-invariant metric $Q_t$ on $S_1^3 \times S_2^3$. But this metric depends on $t$, and hence the fibers are not totally geodesic. This completely describes the metric on $P_{k,l}$. It is interesting to observe that our metrics are just slightly more general than connection-type metrics in that the metrics on the fiber are allowed to depend on a single parameter $t$.

Notice that the left-invariant metric $Q_a$ on G is also right invariant under the maximal torus $T^3_\pm \subset S^3_1 \times S^3_2 \times S^3_3$ containing $K^0_\pm$ which hence acts $S^1$ ineffectively and by isometries on each half $D_\pm = G \times_{K_\pm} D^2$ of $P_{k,l}$. But the intersection of $T_-$ and $T_+$ acts trivially. Furthermore, the first two components of $T^3_\pm$ act by isometries via right translation on the left-invariant metric $Q_t$ on $S^3_1 \times S^3_2$.

We now consider the geometry of the associated bundles $M_{k,l} = P_{k,l} \times_{S^3_1 \times S^3_2} S^3(r)$ and $E_{k,l} = P_{k,l} \times_{S^3_1 \times S^3_2} \mathbb{R}^4$, where we also allow ourselves the freedom of varying the radius in $S^3(r)$. The horizontal distribution and the metric on the base $S^4$ is the same as before, and we only need to describe the metric on the fibers. The fiber of the $S^3$ bundle $M_{k,l} \to S^4$, over the point $c(t) \in S^4$ can be described as $S^3_1 \times S^3_2 \times_{S^3_1 \times S^3_2} S^3(r) = S^3$ where the metric on $S^3_1 \times S^3_2$ is given by the left-invariant metric $Q_t$. Only the right translations on $S^3_1 \times S^3_2$, that are still isometries of $Q_t$, are isometries of this metric on $S^3$. These right translations consist of the action by $T^2 = (e^{i\theta}, e^{i\psi})$ (in the case of $D_-$) and this action of $T^2$ on $S^3$ is the standard cohomogeneity one action on $S^3$. Hence all fibers of $M_{k,l} \to S^4$ are cohomogeneity one metrics on $S^3$ (and in particular not homogeneous). Choosing a basis of $T^2$, the metric on $S^3$ can be described by the length and inner product of the corresponding action fields along a normal geodesic in the fiber. But, unlike in the case of $S^4$, since the principal isotropy group of the $T^2$ action is trivial, the inner product between these two action fields does not have to be 0. Hence the metric on the fiber is described by three functions $h_1(s,t), h_2(s,t), h_3(s,t)$, where $s$ is the arc length parameter of a normal geodesic in the cohomogeneity one metric on the fiber $S^3$ over $c(t)$. Thus, the metric on $M_{k,l}$ is completely described by $(f_1, f_2, f_3, h_1, h_2, h_3): I \times I \to \mathbb{R}^6$.

Similarly, the metric on the fibers of $E_{k,l} \to S^4$ are cohomogeneity two metrics $dt^2 + g_t d\theta^2$ with $g_t$ a cohomogeneity one metric on $S^3$ as above. In both cases, the fibers again change from point to point and hence are not totally geodesic.

Notice that on $M_{k,l}$ (but not on $E_{k,l}$) one can describe a different metric using the identification $P_{k,l} \times_{S^3_1 \times S^3_2} S^3 = P_{k,l}/\triangle S^3$ as a submersed metric from $P_{k,l}$. For this metric, the horizontal distribution and the metric on the base is the same, but the metric on the fibers is now the metric on $S^3 = \triangle S^3 \backslash S^3_1 \times S^3_2$ induced from the left-invariant metric $Q_t$ on $S^3_1 \times S^3_2$ which is invariant under right translations by $T^2$ and hence again only a cohomogeneity one metric on $S^3$. To compare this metric with the previous metric, if we consider the metric on $M_{k,l} = P_{k,l} \times_{S^3_1 \times S^3_2} S^3(r)$ and let r go to infinity, then the limit is the new metric just described. Indeed, if we set $L = S^3_1 \times S^3_2, H = \triangle S^3$, then one has the identification $P \times_L L/H \simeq P/H$ given by $[(p, \ell H)] \to H \cdot \ell^{-1} p$. The fiber $L \times_L L/H$ then gets identified with $H \backslash L$ and if $Q_t(X, Y) = Q(A_t X, Y)$ it follows as in (2.1) that the metric is induced from the left invariant metric $Q(\frac{r^2 A_t}{A_t + r^2} X, Y)$, which as $r$ goes to infinity, converges to $Q_t$.

Next, we consider the isometry group of our examples. As was explained in Section 4, the action of SO(3) on the principal bundle $P_{k,l}$ descends to an action of SO(3) on the total space of $M_{k,l}$ which acts by isometries in the metric of nonnegative curvature that we constructed. Furthermore, the action can be extended to an isometric action by O(3). We suspect that this group will always be the full isometry group of our metrics. In the special case of exotic spheres $M_{k,1-k}$, it follows from [St1, Theorem C], that the id component of the full isometry group can be at most SO(3) × SO(2) or one of its finite quotients. An affirmative answer to the following question would of course rule out such an extension, and would imply that the group SO(3) is always the id-component of the full isometry group.

PROBLEM 5.5. *Does the* SO(3) *subaction of an* (*almost*) *effective* SO(3) × SO(2) *action on an* (*exotic*) 7-*sphere have isotropy groups containing* SO(2)?

In [Da] it was observed that there are natural lifts of the cohomogeneity one action of SO(3)SO(2) on $S^4$ to each of the manifolds $M_{k,l}$. Also, each exotic 7-sphere can be exhibited as a Brieskorn variety, and as such it again supports a natural action of SO(3)SO(2). If the exotic sphere is of the form $M_{k,1-k}$, this action is in general different from the previous ones. In either case, however, the subaction of SO(3) is never almost free.

We also remark, that the action of SO(3)SO(2) on $S^4$ lifts not only to $M_{k,l}$ as in [Da], but to the principal bundle $P_{k,l}$ as well. But this lift does not commute with the free action of $S^3 \times S^3$ and hence one does not obtain a cohomogeneity one action on $P_{k,l}$ as we do in our examples.

An essential difference between our examples and the Gromoll-Meyer metric [GM] on $M_{2,-1}$, is that it has a 4-dimensional isometry group, which agrees with the action of SO(3)SO(2) in [Da]. We finally rephrase the description of this metric on the Gromoll-Meyer sphere in our context, thereby exhibiting similarities and differences.

Consider the following subgroup $G = (S^1 \times S^3_1) \times (S^3_2 \times S^3_3) \subset \mathrm{Sp}(2) \times \mathrm{Sp}(2)$:

$$\left\{ \left( \left( \begin{array}{cc} \cos\theta & -\sin\theta \\ \sin\theta & \cos\theta \end{array} \right) \cdot \mathrm{diag}\,(q_1, q_1), \mathrm{diag}\,(q_2, q_3) \right) \,\Big|\, q_i \in \mathrm{Sp}(1) \right\}.$$

This subgroup acts on Sp(2) (via left and right multiplication) and any two combinations of the $S^3$ factors act freely and one easily sees that in all cases the quotient is $S^4$ on which $S^1 \times S^3$ (the $S^3$ being the remaining $S^3$ factor) acts by cohomogeneity one with the standard sum action. Hence $G$ also acts by cohomogeneity one on $\mathrm{Sp}(2)$ with singular orbits of codimension two and three, and isometrically with respect to the bi-invariant metric with nonnegative curvature. If one chooses the free action by $S_2^3 \times S_3^3$, the principal $S^3 \times S^3$ bundle over $S^4$ is the two-fold cover of the frame bundle of the tangent bundle of $S^4$ and hence $P_{1,1} = \mathrm{Sp}(2)$. If one chooses the free action by $S_1^3 \times S_2^3$, it was shown in [GM] that one obtains the principal bundle $P_{2,-1}$ and hence the associated sphere bundle $P_{2,-1} \times_{S_1^3 \times S_2^3} S^3 = P_{2,-1}/\triangle_{1,2}S^3$ is the exotic sphere $M_{2,-1}$ with a submersed metric of nonnegative curvature. As in our case, the action of $S^1 \times S_3^3$ descends to an action of the associated $S^3$ bundle $M_{2,-1}$ which becomes an effective action by $\mathrm{SO}(3)\mathrm{SO}(2)$ and which, by [St1, Theorem C], is the id-component of the isometry group of the Gromoll-Meyer metric. Notice that, as in our case, we also get a family of metrics with nonnegative curvature on the Gromoll-Meyer sphere by considering $\mathrm{Sp}(2) \times_{S_1^3 \times S_2^3} S^3(r)$, and as $r$ goes to infinity we obtain the Gromoll-Meyer metric in the limit.

Of course, as a consequence of our results, it follows that $P_{2,-1} = \mathrm{Sp}(2)$ also has infinitely many cohomogeneity one actions by $S^3 \times S^3 \times S^3$, but with singular orbits of codimension two. Another major difference between our metrics induced by these actions and the Gromoll-Meyer metric, is that in their example there exists an open set of points in $M_{2,-1}$ on which every two-plane has positive curvature, whereas in our example, by construction, there are always 2-planes of $0$ curvature at every point.

Motivated by Proposition 4.8 and Theorem G we conclude our discussion with the following natural question.

PROBLEM 5.6.    *On each of the Milnor spheres (including the standard sphere), as well as on the homotopy $\mathbb{R}P^5$ in Theorem G, does the space of metrics with nonnegative sectional curvature have infinitely many components?*

In a similar vein, in [KS] it was shown that on some of the homogeneous spaces $\mathrm{SU}(3)/S^1$ the space of metrics with positive sectional curvature has at least two components.

*Added in proof.* Motivated by Problem 5.2, D. Crowley and C. Escher [CE] recently completed the classification of the total spaces $M_{k,l}$ of $S^3$ bundles over $S^4$ up to orientation preserving and reversing homotopy equivalence, homeomorphism and diffeomorphism.

In [KiS] N. Kichloo and K. Shankar partially answered another question raised in Section 5 by showing that the positively curved Berger space $B^7 = \mathrm{SO}(5)/\mathrm{SO}(3)$ is PL-homeomorphic to an $S^3$ bundle over $S^4$. It remains to compute the Eells-Kuiper invariant for this space in order to see if it is diffeomorphic to one.

It also follows from [CE] that all $M_{k,l}$ with $k + l = 10$ are homotopy equivalent to each other. Hence for the total spaces $P_k$ of principal $S^3$ bundles over $S^4$, it follows that $P_{10}$ is homotopy equivalent to the homogeneous space $B^7$, whereas, as we observed in Section 5, all other $P_k$, $k \geq 3$ are strongly inhomogeneous.

University of Maryland, College Park, MD
*E-mail address:* kng@math.umd.edu

University of Pennsylvania, Philadelphia, PA
*E-mail address:* wziller@math.upenn.edu